%% file: Combined_Wrapper.tex
\documentclass{imanum}
\jno{drnxxx}
\received{21 October 2010}

\usepackage{lmodern,enumerate}

\include{commonheader}
\usepackage{array}
\usepackage{multirow}

\newcommand{\cB}[1]{\mathcal{B}\left[ #1 \right]}
\newcommand{\cc}[2]{#2}
\defcitealias{companion1}{\textit{ibid}}
\begin{document}

\title{High-order finite elements on pyramids:\\approximation spaces, unisolvency and exactness.}

\author{{\sc Nilima Nigam \thanks{nigam@math.sfu.ca. The work of NN was supported by the Natural Sciences and Engineering Research Council of Canada, and the Canada Research Chairs program}} \\[2pt]Dept. of Mathematics, Simon Fraser University, BC, Canada\\[6pt]
{\sc Joel Phillips \thanks{joel.phillips@reading.ac.uk.  JP was supported by a Natural Sciences and Engineering Research Council graduate fellowship.}} \\[2pt]Dept. of Mathematics and Statistics, Reading University, UK. }

\maketitle

 
\begin{abstract}
{We present a family of  high-order finite element approximation spaces on a pyramid, and associated unisolvent degrees of freedom. These spaces consist of rational basis functions. We establish  conforming, exactness and polynomial approximation properties.}
{conforming finite elements; pyramid; high order approximation, commuting diagram}
\end{abstract}

\include{highorderpyramid-part1}
\include{highorderpyramid-part2}

{\bf Acknowledgement}
We gratefully acknowledge the contributions of Leszek Dem\-kowicz, who suggested this problem.  We would like to thank Leszek Demkowicz,  Peter Monk  and Paul Tupper for helpful discussions on the paper.  The work of Nilima Nigam was partially supported by NSERC,  FQRNT, and the Canada Research Chairs program. Joel Phillips was supported by NSERC.

\appendix
\include{Combined-tables-appendix}

\bibliographystyle{IMANUM-BIB}

\bibliography{pyramid}

\end{document}

%% file: commonheader.tex
\usepackage{amsmath,enumerate}
\usepackage{amsfonts}
\usepackage{amssymb}
\usepackage{epsfig}

\usepackage{graphicx}
\usepackage{amscd}

\newcommand{\x}{\mathbf x}
\newcommand{\xx}{x}

\newcommand{\xyy}{y}
\newcommand{\xz}{z}
\newcommand{\xix}{\xi}
\newcommand{\xiy}{\eta}
\newcommand{\xiz}{\zeta}

\newcommand{\ba}{\begin{align} }
\newcommand{\bs}{\begin{split} }
\newcommand{\es}{\end{split} }
\newcommand{\ea}{\end{align} }

\newcommand{\xii}{\boldsymbol{\xi}}

\newcommand{\uu}{u}
\newcommand{\vv}{v}

\newcommand{\G}{\mathcal{U}^{(0),k}}
\newcommand{\C}{\mathcal{U}^{(1),k}}
\newcommand{\D}{\mathcal{U}^{(2),k}}
\newcommand{\Z}{\mathcal{U}^{(3),k}}
\newcommand{\HH}{H}

\newcommand{\PPP}{P}
\newcommand{\RR}{\mathbb R}
\newcommand{\ip}{\Omega_{\infty} }
\providecommand{\abs}[1]{\lvert#1\rvert}
\providecommand{\col}[3]{\begin{pmatrix} #1 \\ #2 \\ #3 \end{pmatrix} }
\providecommand{\coll}[2]{\begin{pmatrix} #1 \\ #2 \end{pmatrix} }

\providecommand{\pd}[2]{\frac{\partial #1}{\partial #2}}

\providecommand{\set}[1]{\left\{ #1 \right \}}

\DeclareMathOperator{\curl}{{curl}}
\DeclareMathOperator{\spnn}{{span}}
\DeclareMathOperator{\grad}{{grad}}
\DeclareMathOperator{\divv}{div}

\providecommand{\spn}[1]{\spnn \left\{ #1 \right\} }

\usepackage{graphicx}

%% file: highorderpyramid-part1.tex
\section{Introduction}

High order conforming finite elements for $H(\curl)$ and $H(\divv)$ spaces based on meshes composed of tetrahedra and hexahedra were first presented by \cc{Nedelec, \cite{nedelec:mixed1980}}{\citet{nedelec:mixed1980}}. The demands of the specific problem geometry (regions with complex features as inclusions) or efficient calculation (design of unstructured hexahedral meshes) may necessitate the use of hybrid meshes which include both tetrahedral and hexahedral elements, see e.g. \cite{bergotcohendurufle::pyramid}.   If these meshes are to avoid hanging nodes then they will, in general, contain pyramids.  A hybrid mesh may contain tetrahedra to provide
localised h-refinement and computationally-efficient cuboids to fill large spaces in which the solution is regular, and pyramids to glue these together. This situation is nicely illustrated in \cite{owen:hextet}. Note that (triangularly) prismatic elements are also required; these turn out to be
relatively straightforward to construct and analyse, see e.g. \cite{chen1989prismatic}. Pyramidal elements also arise more explicitly when attempting to mesh thin three dimensional structures using prismatic elements, see \cite{lee2005thinwalled,gatto2010construction}.

Consider a contractible domain $D\in \mathbb{R}^3$ which is triangulated using a mesh containing both tetrahedral and hexahedral elements. If one is to avoid hanging nodes or edges, the triangulation must also, in general, include quadrilateral-based pyramids.   In what follows, we assume these pyramids can be mapped in an affine manner to a reference pyramid, $\Omega$, which has  a square base and is defined as:
\begin{align}\label{pyramid_def}
 \Omega = \lbrace \xii=   (\xi,\eta,\zeta) \in \RR^3  \; | \; \xi,\eta,\zeta \geq 0, \; \xi \leq 1-\zeta, \; \eta \leq 1-\zeta\rbrace. 
\end{align} 
It is our aim to construct high order finite elements on such a pyramid. Concretely, in this paper we present finite element triples, $(\Omega, \mathcal{U}^{(s),k}(\Omega), \Sigma^{(s),k} )$, for positive integers $k$ which are unisolvent conforming finite elements for  $H^1(\Omega)$, $H(\curl,\Omega)$, $H(\divv,\Omega)$ and $L^2(\Omega)$ respectively for $s=0,1,2,3$.  Here $\mathcal{U}^{(s),k} (\Omega)$ denotes the $k$th order finite dimensional approximation space for the relevant Sobolev space and the sets $\Sigma^{(s),k}$ are the associated degrees of freedom.  We seek finite elements with the following properties {\bf P1-P3}:
\begin{enumerate}[{\bf P1)}]
\item \emph{Compatibility:} Not only should the elements be conforming, but the restriction of each element to its triangular and quadrilateral face(s) should match that of the corresponding canonical tetrahedral and hexahedral finite element.  This means that both the spaces spanned by the traces and the external degrees of freedom on faces and edges are the same as those of the usual tetrahedral/hexahedral elements (see Table \ref{tab:dof}).  In other words, the elements should satisfy the correct \emph{patching conditions} on inter-element boundaries, \cc{\cite{hiptmair:pyramid}}{\citep[see][]{hiptmair:pyramid}}.  We will use \cite{monk:fe-maxwell} as our reference for the tetrahedral and hexahedral spaces and external degrees of freedom, see Table 1.

\item \emph{Approximation:} The discrete spaces $\mathcal{U}^{(s),k}(\Omega)$ should allow for high-order approximation to the spaces $H^1(\Omega), \HH(\curl,\Omega)$, etc. In particular, given a positive integer $p$, it should be possible to choose $k$ such that all polynomials of degree $p$ (we denote these by $P^k\equiv P^k(\Omega)$) are contained in $\mathcal{U}^{(s),k}(\Omega)$.

\item \emph{Stability:} The elements satisfy a commuting diagram property:
\small
\begin{equation} \label{cd}
\begin{CD}
H^r(\Omega) @>\nabla>> \HH^{r-1}(\curl, \Omega) @>\nabla\times>> \HH^{r-1}(\divv, \Omega) @>\nabla\cdot>> H^{r-1}(\Omega) \\
@V \Pi^{(0)} VV @V \Pi^{(1)} VV @V \Pi^{(2)} VV @V \Pi^{(3)} VV \\
\G (\Omega) @>\nabla>> \C (\Omega) @>\nabla\times>> \D (\Omega) @>\nabla\cdot>> \Z(\Omega)
\end{CD}
\end{equation}
\normalsize
 Here $\Pi^{(s)}$, s=0,1,2,3, denote interpolation operators induced by the degrees of freedom, $\Sigma^{(s),k}$ and $r$ is chosen so that the interpolation operators are well defined.
\end{enumerate}

 \begin{table}[htbp]
   \centering
   \begin{tabular}{| c|>{\small}c|  >{\small}c| >{\small}c |}  \hline

  & {Edge $e$ } & \multicolumn{2}{c|}{Face $f$ }  \\ \cline{2-4}
   &  tetrahedra  &tetrahedra & hexahedra \\ 
  & \& hexahedra & & \\\hline
   
    $H^1(\Omega)$ & $ \int_e pq \,ds $ &  $ \int_f pq \, dA$& $ \int_f pq \, dA$   \\ 
       
   &  $\forall q\in P^{k-2}(e)$  &$\forall q \in P_{k-3}(f)$& $ \forall q \in Q^{k-2,k-2}(f)$ \\  \hline

     $  \HH(\curl,\Omega)$   &$\int_e{\bf u}\cdot {\mathbf t} q\, ds$  &$\int_f {\bf u}\cdot{\bf q}\,dA$ & $\int_f {\bf u}\times \nu \cdot{\bf q}\,dA$\\
      
      & $\forall q \in P^{k-1}(e)$&  $\forall {\bf q} \in P^{k-2}(f), {\bf q} \cdot \nu =0$& ${\bf q} \in Q^{k-2,k-1}\times Q^{k-1,k-2}(f)\,\, $\\  \hline
     
     $ \HH(\divv,\Omega)$  & \---&$\int_f {\bf u}\cdot {\mathbf \nu} q \,dA$ & $\int_f {\bf u}\cdot {\mathbf \nu} q \,dA$ \\
       & & $\forall q \in P^{k-1}(f)$& $\forall q \in Q^{k-1,k-1}(f)$ \\
      \hline
       
         \hline
   \end{tabular}
   \caption{Edge and face degrees of freedom for tetrahedral and hexahedral reference elements. The vertex degrees of freedom for the $H^1(\Omega)$ elements on tetrahedra and hexahedra are the same. There are no exterior degrees of freedom for the $L^2(\Omega)$ approximants. $t$ is the unit tangent along an edge, and $\nu$ the unit outer normal to a face. {\it We denote by $P^k$ polynomials of maximal degree $k$; we employ the same notation for scalar and vector-valued objects.} }
   \label{tab:dof}
\end{table}



\cc{Gradinaru and Hiptmair }{}\cite{hiptmair:pyramid} constructed ''Whitney'' elements satisfying properties {\bf P1} and {\bf P3} and our family of elements includes these as the lowest order case, see Section \ref{firstorder}.  In the engineering literature, \cc{Zgainski et al }{}\cite{zgainski:edge,zgainski:family} appear to have discovered the same first order $H(\curl)$-conforming element independently and also demonstrated a second order element. \cc{In }{}\cite{bergotcohendurufle::pyramid}\cc{ the authors}{} describe high-order finite elements for $H^1(\Omega)$, but not the other spaces. \cc{Graglia et al }{}\cite{graglia:highorderpyramid} constructed $H(\curl)$ and $H(\divv)$ elements of arbitrarily high order.   Similarly, \cc{Sherwin }{}\cite{sherwin:hierarchicalhp} demonstrated $H^1$-conforming elements also satisfying properties (1) and (2).  These high order constructions provide an explicit scheme for determining nodal basis functions; none of them address the commuting diagram property, {\bf P3}.  

The mimetic finite difference method, originally presented in \cc{1997 }{}\cite{hyman:naturaldiscretizations} and further developed by several authors \cc{(e.g., \cite{kuznetsov:mimetic,campbell:mimetic,brezzi:familymimetic,brezzi:mimeticconvergence})}{\citep[e.g.][]{kuznetsov:mimetic,campbell:mimetic,brezzi:familymimetic,brezzi:mimeticconvergence}}  develops low-order approximations on polyhedral meshes and hence includes pyramids as a special case.

Our starting point is an observation: that it is not always possible to extend polynomial data on the faces of a pyramid using a polynomial within the pyramid.  Indeed, it is {\it impossible} to construct useful $H^1(\Omega)$ pyramidal finite elements using only polynomial basis functions. Specifically, in Theorem~\ref{impossible}, we demonstrate an $H^1(\Omega)$ function which has polynomial traces on the faces of the pyramid, but which does not admit a polynomial representation in the pyramid itself.

\begin{theorem}\label{impossible} Let $\Omega$ be the pyramid defined in \eqref{pyramid_def}. Consider the function $u:\Omega\rightarrow \RR$ defined by $$u(\xix,\xiy,\xiz) = \frac{\xix\xiz(\xix+\xiz-1)(\xiy+\xiz-1)}{1-\xiz}.$$  Then, \begin{enumerate}
\item $u \in H^1(\Omega)$,
\item  $u$ has polynomial traces on the pyramid faces,
\item  $u$ cannot be represented by any polynomial function on $\Omega$ which also satisfies property ({\bf P1}).
\end{enumerate}
\end{theorem} 
\begin{proof} 
It is straightforward to verify (1). It is easy to see $u\vert_{\xiy=0}= -\xix \xiz(\xix+\xiz-1)$ and $u = 0$ on the other faces of the pyramid. This establishes (2).

We prove (3) by contradiction.  Since $\Omega$ has Lipschitz boundary, we can extend $u$ to a function $U \in H^1(\RR^3)$ \cc{(see, for example, \cite{adams:sobolev})}{\citep[see, for example][]{adams:sobolev}}.  Suppose that we could represent $u=U\vert_\Omega$ by a polynomial function $p$, in a manner consistent with property ({\bf P1}).   The traces of $U$ on the faces will then be interpolated exactly by the polynomial Whitney forms specified by adjacent neighbouring tetrahedra and hexahedra. Since an $H^1$-conforming approximation must be continuous across interelement faces, we must have $p=U$ on each face of the pyramid.

Since $U=u=0$ on four of the faces of the pyramid, we can factorise: 
\begin{align}
p(\xix,\xiy,\xiz)&= \xix\xiz(\xix+\xiz-1)(\xiy+\xiz-1)\left[r(\xix,\xiz)+\xiy s(\xix,\xiy,\xiz)\right], \label{poly1}
\end{align}
where $r$ and $s$ are polynomial.
Further,  $U= -\xix \xiz(\xix+\xiz-1)$ on the face $\xiy=0$ and so: 
\begin{align}
p(\xix,0,\xiz)= \xix\xiz(\xix+\xiz-1) (\xiz-1) r(\xix,\xiz) = -\xix \xiz(\xix + \xiz -1) ,\label{poly2}
\end{align}
which implies  that $(\xiz-1)r(\xix,\xiz) = -1$. This contradicts the polynomial nature of $r$.\end{proof}
 
A similar result is presented \cc{in \cite{wieners:conforming}}{by \citeauthor{wieners:conforming}}, where it is claimed that, under the assumption that shape functions must be polynomial, there exists no continuously differentiable conforming shape functions for the pyramid which are linear / bilinear on the faces.

The insufficiency of polynomials can be seen in all previous successful attempts to construct pyramidal finite elements.  In addition to  \cite{hiptmair:pyramid}, finite element bases that include rational functions are given by, e.g., \cite{graglia:highorderpyramid,sherwin:hierarchicalhp,zgainski:edge,zgainski:family} and \cc{ \cite{felippa:compendium,owen:hextet,wieners:conforming, davies:symmetric}}{ \citeauthor{wieners:conforming}. In \cite{felippa:compendium,owen:hextet,davies:symmetric},} the authors use piecewise polynomial functions via a macro-element that divides the pyramid into two or four tetrahedra.  Interestingly, although \cite{wachspress} only applies his construction to a class of polyhedra that does not include pyramids, this restriction appears to be unnecessary and the ``rational finite elements'' given therein appear to include the high order $H^1$ pyramidal elements as a special case.

The major contribution of this paper is a comprehensive development of high-order finite elements on a pyramidal element. We will present candidate approximation spaces $\mathcal{U}^{(s),k}(\Omega)$ for $s=0,1,2,3$ and $k \in \mathbb{N}$, by first developing these on an infinite reference pyramid.  We also show that these spaces admit convenient Helmholtz-like decompositions, and that their traces on faces and edges are consistent with traces from neighbouring elements. Hence property {\bf P1} is satisfied by $\mathcal{U}^{(s),k}(\Omega)$. As a concrete example, we verify that our first order elements agree with those presented in \cite{hiptmair:pyramid}.

Next, we provide a description of the degrees of freedom, $\Sigma^{(s),k}$ and demonstrate unisolvency.  The exterior degrees of freedom agree precisely with those specified by neighbouring tetrahedral or hexahedral elements.  Properties {\bf P2} and {\bf P3} are also established. We will use the projection-based interpolation described in \cite{demkowicz:projection,demkowicz:derham} to solve the difficult problem of defining the internal degrees of freedom on a pyramid.  It is possible to use projection based interpolation for the external degrees too, and we believe that the $hp$ framework of which it is a part will also accommodate our element.  However, this is not our immediate objective and the external degrees described in \cite{monk:fe-maxwell} allow for a more explicit exposition.

In Section 6, we show that the discrete spaces  $\mathcal{U}^{(s),k}(\Omega)$ form an exact subcomplex. That is, we show that $d \mathcal{U}^{(s),k}(\Omega) \subset  \mathcal{U}^{(s+1),k}(\Omega)$ for $s=0,1,2$, and that any discrete $(s+1)$\---form which has a vanishing exterior derivative is derivable from a discrete potential which is an $s$\---form. 
These spaces,  along with the interpolants which are induced by the degrees of freedom, satisfy a ``commuting diagram property'' which is crucial to the stable computation of mixed problems. Finally we show that these finite elements are indeed high-order in the sense that they include high-degree polynomials. While the inclusion of high-degree polynomials is an important step towards approximability, we shall show in a subsequent paper that the usual finite element arguments need modification in our context. In particular, since the spaces $\mathcal{U}^{(s),k}(\Omega)$ contain rational functions, it is not true that high derivatives evaluate to 0, in sharp contrast to the situation for polynomials.

The organization of the rest of this paper is as follows:
\begin{enumerate}[{\bf Section} 1] \setlength{\itemsep}{0cm}
\setcounter{enumi}{1}
\item{\it The infinite reference element: some preliminaries}
\item{\it The approximation spaces  $\mathcal{U}^{(s),k}(\ip)$ on the infinite pyramid}
\item{\it The approximation spaces  $\mathcal{U}^{(s),k}(\Omega)$ on the finite pyramid}
\item {\it The degrees of freedom $\Sigma^{(s),k}$ and unisolvency}\label{section:unisolvency}
\item {\it Interpolation and exact sequence property}
\item {\it Polynomial approximation property}
\item {\it Appendix: Shape functions}
\end{enumerate}

\section{The infinite reference element: some preliminaries}
As discussed earlier, pyramidal finite element spaces must include rational functions.  
To construct the finite elements, we shall make use of two reference elements: the {\it finite pyramid}, $\Omega$, already introduced in \eqref{pyramid_def}, and the {\it infinite pyramid} $\ip$.  The infinite pyramid is an unusual, but instructive domain; it possesses hexahedral symmetries which will allow us to specify or study important properties for the approximation spaces. We will then map these spaces to the finite pyramid.

 We will typically use the symbols $(\xx,\xyy,\xz)$ as coordinates on the infinite pyramid  and $(\xix,\xiy,\xiz)$ on the finite pyramid. 
The infinite reference pyramid is defined as 
\begin{align}
\ip = \lbrace \x = (\xx,\xyy,\xz) \in \RR^3 \cup \infty \; | \; \xx,\xyy,\xz \geq 0, \; \xx \leq 1, \; \xyy \leq 1 \rbrace. \label{infinite_pyramid}
\end{align}

\begin{figure}
  \includegraphics[width=2.2in]{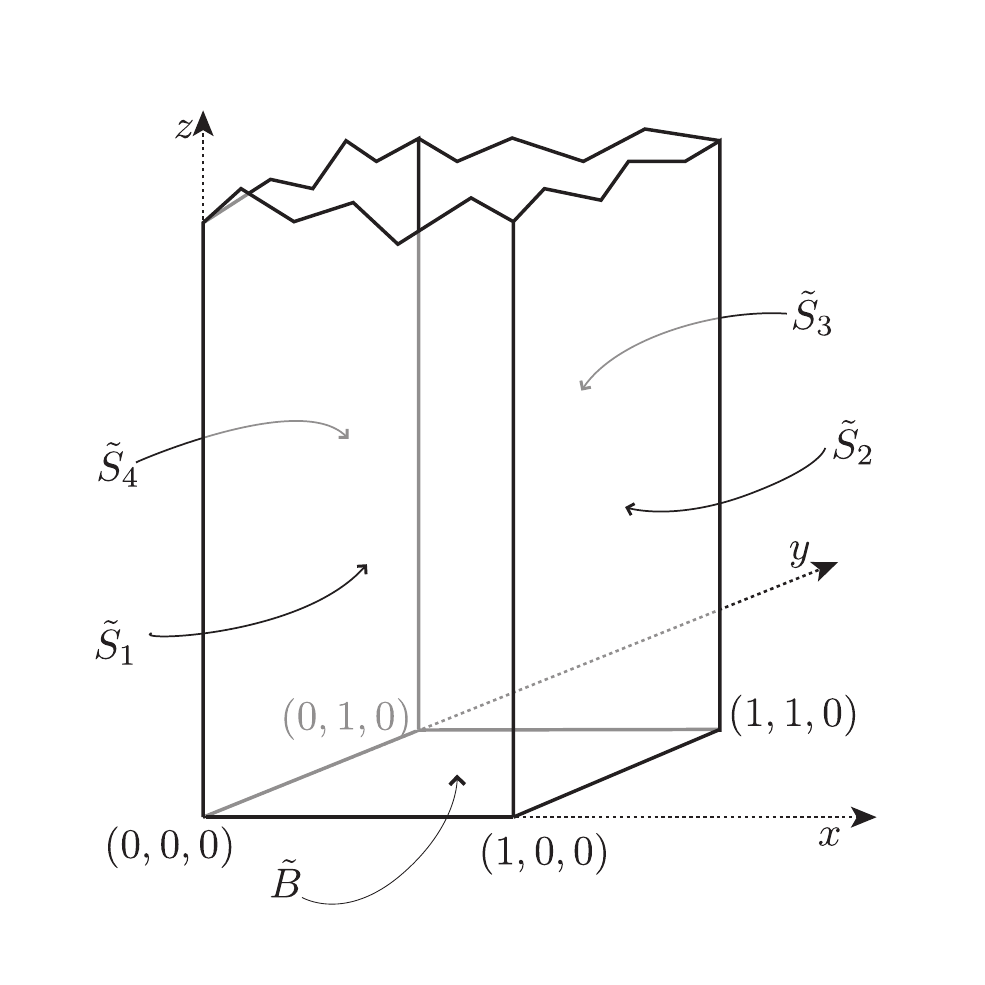}\qquad\qquad \includegraphics[width=2.2in]{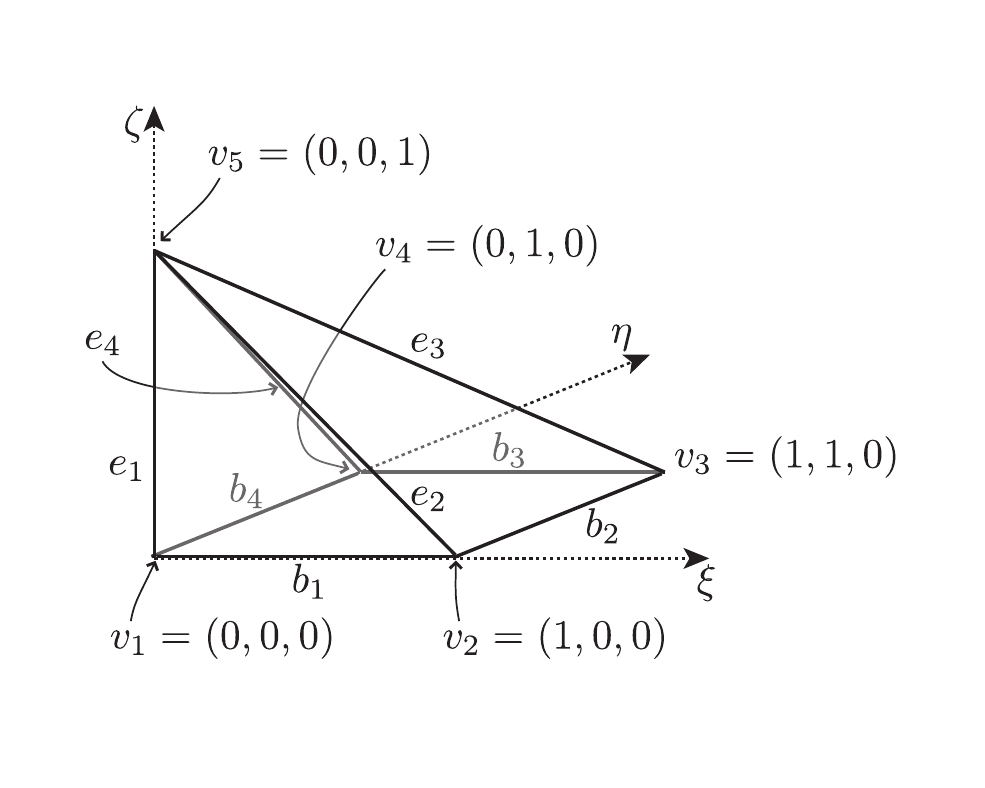}

\caption{Left: The infinite pyramid $\ip$. Right: The finite reference pyramid $\Omega$}
\label{fig:infpyramid}
\end{figure} 

Figure \ref{fig:infpyramid}  shows the two pyramids.  The vertical faces of the infinite pyramid lie in the planes $y=0$, $x=1$, $y=1$, $x=0$.  We denote them as $S_{1,\ip}, S_{2,\ip}, S_{3,\ip}$, and $S_{4,\ip}$ respectively, and the corresponding faces on the finite pyramid $S_{i,\Omega} = \phi(S_{i,\ip})$.  Let $B_{\ip}$ refer to the base face, $z=0$, of the infinite pyramid and $B_\Omega$ the base face of the finite pyramid.  The vertices of the finite pyramid are denoted $v_i$, $i=1..5$, with $v_5$ the point $(0,0,1)$.  
\subsection{The infinite reference element: pullbacks}\label{infiniteref}
To associate the finite and infinite pyramids, define the bijection  $\phi: \ip \rightarrow \Omega$
\begin{equation}\label{phi}
\begin{split}
&\phi(\xx,\xyy,\xz) = \left(\frac{\xx}{1+\xz}, \frac{\xyy}{1+\xz}, \frac{\xz}{1+\xz} \right), \qquad \phi(\infty) = (0,0,1),
\end{split}
\end{equation}
which is a diffeomorphism if we restrict the domain to $\ip \backslash \infty$ (and the range to the finite pyramid with its tip removed).


The infinite pyramid will serve as a  tool for the construction of the function spaces for the elements.  We thus need to understand how to map functions between spaces on the finite pyramid, $\mathcal{U}^{(s),k}(\Omega)$ and the infinite pyramid, $\mathcal{U}^{(s),k}(\ip)$.  A major consideration is for approximation spaces on the infinite pyramid to satisfy an exact sequence property.  To have this exact sequence property preserved on the finite pyramid, it is necessary that the mappings between the spaces on the finite and infinite pyramids commute with the $\grad$, $\curl$ and $\divv$ operators.

In the language of differential geometry, where the elements of each space can be considered to be proxies for 0, 1, 2 and 3-forms, the mappings should be pullbacks. We shall use the same notation for each map - the context will never be ambiguous.  We point the reader to \cite{arnold:acta} for an excellent treatment of the finite element exterior calculus.  In this paper, we will switch between referring to objects as forms or functions, depending on the context.  Formally (because we have not yet defined the appropriate Sobolev spaces on the infinite pyramid):
\begin{subequations}
\begin{align}
&\forall u \in H^1(\Omega) & \phi^*u  &= u \circ \phi, \\
&\forall E \in \HH(\curl, \Omega) & \phi^*E  &= D\phi^T \cdot[E \circ \phi], \\
&\forall \vv \in \HH(\divv, \Omega) & \phi^*\vv  &= \abs{D\phi}D\phi^{-1} \cdot [\vv \circ\phi] ,\\
&\forall q \in L^2(\Omega) & \phi^*q  &= \abs{D\phi}[q\circ\phi] ,
\end{align}
\end{subequations}
where $D\phi$ is the Jacobian matrix, $\frac{1}{(\xz+1)^2}\col{\xz+1 & 0 & -\xx}{0 & \xz+1 & -\xyy}{0 & 0 & 1}$.  The pullback is a bijection and the inverse pullback, $(\phi^*)^{-1}$ is equal to $(\phi^{-1})^*$. Since $z\geq 0$, $D\phi^T D\phi$ is positive definite.

The infinite reference element is a convenient tool, since it possesses both rotational symmetries and the tensorial nature of regular hexahedral elements. This is particularly useful while discussing traces onto the boundaries of the pyramid.

\subsection{The infinite reference element:  Sobolev spaces}
The infinite reference pyramid has obvious symmetries, which make it easier to specify and analyze approximation spaces. However, it has semi-infinite extent along the $z$-direction, and we must therefore describe analogues of $H^1(\Omega)$, $H(\curl, \Omega)$ etc. on $\ip$. Not surprisingly, these Sobolev spaces will have weighted norms.
\begin{definition} \label{infinitesobolev}
Let $\Omega_\infty$ be the infinite pyramid defined in \eqref{infinite_pyramid},  and   $\phi: \ip \rightarrow \Omega$ be the pullback map. We define the following inner product spaces:

 $H^1_w(\Omega_\infty)$ is the closure of the  set of smooth scalar-valued functions $v:\ip\rightarrow \RR$ under the norm induced by the inner product$$ (u,v) _{H_w^1(\ip)} := \int_{\ip} \frac{u v}{(1+z)^4} + (\nabla u)^T \mathcal{A} \nabla v d \x. $$ Here $\mathcal{A} = |D\phi| D\phi^{-1}{D\phi^{-1}}^T$ is positive definite. 
$H_w(\curl, \Omega_\infty)$ is  the closure of the  set of smooth vector-valued functions (1-forms) $F:\ip\rightarrow (\RR)^3$ under the norm induced by  inner product 
$$ (F,G) _{H_w(\curl,\ip)} := \int_{\ip}   (F)^T \mathcal{A} (G) + (\curl F)^T\mathcal{B} (\curl G) d \x. $$ Here $\mathcal{B} = |D\phi^{-1}|D\phi^T D\phi$,  and is positive definite.  $H_w(\divv, \ip)$ is the  closure of the set of smooth vector-valued functions (2-forms) $F:\ip\rightarrow (\RR)^3$ with inner product 
$$  (F,G)_{H_w(\divv,\ip)} := \int_{\ip}   (F)^T \mathcal{B} (G) + (\divv F)^T(1+z)^4(\divv G) d \x. $$
$L_w^2(\ip)$ is the  closure of the set of smooth scalar-valued functions (3-forms) with inner product,
$$  (u,v)_{L_w^2(\ip)} := \int_{\ip}    (1+z)^4(uv) d \x. $$

\end{definition}

 \begin{remark}
 We  observe that the inner products  on the infinite pyramid are weighted by powers of $\frac{1}{(1+z)}$. The subscript $w$ is used to emphasize that these are weighted norms. The weights are entirely specified by the projective mapping, $\phi$, and the associated pull-backs for the various forms. 
It is important to note, for example, that 
$ \|u\|^2_{L_w^2(\ip)} = \int_{\ip} \frac{u^2}{(1+z)^4} d\x$ if $u$ is a zero form, while $ \|u\|^2_{L_w^2(\ip)} = \int_{\ip} {u^2}{(1+z)^4} d\x$ if $u$ is a 3-form. 

\end{remark}
These inner product spaces can be related to more familiar Sobolev spaces on the finite pyramid, as is done in the following theorem: 
\begin{lemma}\label{identification_theorem} It is easy to verify that 
the inner product spaces $H_w^1(\ip), H_w(\curl,\ip), H_w(\divv,\ip)$ and $L_w^2(\ip)$ in Definition \ref{infinitesobolev} are Hilbert spaces. Morever, $\phi^* : H^1(\Omega) \rightarrow H_w^1(\ip)$ is an isometry.
The analogous statements are true for $H_w(\curl,\ip)$,$H_w(\divv, \ip)$ and $L^2_w(\ip)$.

\end{lemma}
\begin{proof}
The pullbacks, $\phi^*$ are formally bijections because $\Omega$ and $\ip$ have the same dimension.  Suppose $\tilde{u} $ is a 0-form in $H^1(\Omega)$ and let $u = \phi^*\tilde u$.  Then
$$ \|\tilde u\|^2_{L^2(\Omega)} =  \int_{\ip} |D\phi| |{u}(\x)|^2 d \x = \int_{\ip} \frac{1}{(1+\xz)^4} |{u}(\x)|^2 d\x. 
$$
Now, the gradient and pull-back operator commute. We can thus use the appropriate pull-back to obtain
\begin{align*}
\|\nabla \tilde{u}\|_{L^2(\Omega)}^2 &= \int_\Omega |\nabla \tilde{u}|^2 d \xii 
=\int_\Omega |{D\phi^{-1}}^T \nabla {u} \circ \phi^{-1}|^2 d \xii \\
&=\int_{\ip} |D\phi| |{D\phi^{-1}}^T \nabla {u}|^2 d \x = \int_{\ip} \nabla{u}^T\mathcal{A}\nabla u d\x. 
\end{align*} 
Hence $\| \tilde{u}\|_{H^1(\Omega)}^2 = \| \tilde{u}\|_{L^2(\Omega)}^2 + \|\nabla \tilde{u}\|_{L^2(\Omega)}^2  = \|u\|_{H^1_w(\ip)} $.  The proofs for $H_w(\curl,\ip)$, $H_w(\divv,\ip)$ and $L_w^2(\ip)$ follow analogously.
\end{proof}

We collect here, for convenience, concrete instantiations of the inverse pullback mapping. 
\begin{subequations}
\begin{align}
&\forall u \in H_w^1(\Omega_\infty), & (\phi^{-1})^*u  &= u \circ \phi^{-1}, \\
&\forall E \in H_w\curl, \ip), & (\phi^{-1})^*E  &= [(1+z)\col{1 & 0 & 0}{ 0&1&0}{ \xx & \xyy & 1+\xz} \cdot E ] \circ \phi^{-1}, \\
&\forall \vv \in H_w(\divv, \ip), & (\phi^{-1})^*\vv &= [(1+\xz)^2 \col{1+z&0&-x}{0&1+z&-y}{0&0&1}\cdot \vv] \circ\phi^{-1},\\
&\forall q \in L_w^2(\ip) ,& (\phi^{-1})^*q  &= [(1+\xz)^4 q]\circ\phi^{-1}.
\end{align}
\end{subequations}
\subsection{Rotations and traces} \label{section:rotations}
Define $R_{\ip}:\ip \rightarrow \ip$ to be the affine mapping that sends the infinite pyramid to itself and rotates it a quarter turn about the axis $x=y=\frac{1}{2}$, that is, the vertical face ${S}_{1,\ip}$ is mapped to ${S}_{2,\ip}$,  the face ${S}_{2,\ip}$ is mapped to ${S}_{3,\ip}$, etc. Explicitly, 
\begin{equation}\label{rotation}
R_{\ip}:(x,y,z) \mapsto (1-y,x,z). 
\end{equation}
We can also define a mapping that sends the finite pyramid to itself, rotating the faces, $R:\Omega \rightarrow \Omega$ by 
\begin{equation*}
R=\phi \circ R_{\ip} \circ \phi^{-1}, \qquad R:(\xix,\xiy,\xiz) \mapsto (1-\xiy-\xiz, \xix, \xiz).
\end{equation*}
It is clear that if an approximation space $\mathcal{U}^{(s),k}(\ip)$ is invariant under the mapping $R_{\ip}$, its (inverse) pullback to the finite pyramid will be invariant under $R$.  This property will prove convenient when we consider exterior shape functions and exterior degrees of freedom. 

The trace map from a manifold to a submanifold is the pullback of the inclusion map for differential forms \cc{(see, for example, \cite{arnold2009}, pg 41 ff.)}{\citep[see, for example,][pg 41 ff.]{arnold2009}} and so we expect that zero trace data will be preserved by the pullback mapping.  The following lemma makes this explicit in our concrete vector calculus formulation, where traces for 1-forms consist only of the tangential components and for 2-forms the normal components.  We suppose that $S_{\ip}$ is a surface of the infinite pyramid and let $S_{\Omega}$ be its image under $\phi$ on the finite pyramid.

\begin{lemma}\label{tracemaps}
\begin{itemize} 
\item A  1-form $u$ is normal to $S_\Omega$ at a point $\xi = \phi(x)$ if and only if the pullback  $\phi^* u$ is normal to $S_{\ip}$ at $x$.
\item A 2-form $u$ is tangent to $S_\Omega$ at a point $\xi=\phi(x)$ if and only if the pullback  $\phi^*u$ is tangent to $S_{\ip}$ at $x$.  
\end{itemize}
\end{lemma}
\begin{proof}
Let $S_{\Omega}$ be described (locally) by $S_\Omega = \lbrace \xi : f(\xi) = 0\rbrace$.  Define $g = f\circ\phi$, then $S_{\ip} = \lbrace x : g(x) = 0\rbrace$. 
To establish the first result,  let $u$ be a 1-form which is normal to $S_{\Omega}$ at $\xi$, then 
\begin{align}\label{ulambdaf}
u(\xi) = \lambda(\xi) \nabla f (\xi)
\end{align}
for some scalar function $\lambda$.  By the chain rule, and substituting \eqref{ulambdaf}
\begin{align*}
\nabla g(x) = (D\phi)^T(x)\cdot (\nabla f)(\phi(x))=(D\phi)^T(x)\cdot \frac{u(\phi(x))}{\lambda(\phi(x))} &= \frac{\phi^* u}{\lambda(\phi(x))} 
\\\Rightarrow\lambda(\phi(x))\nabla g(x) &= \phi^*u(x).
\end{align*}
Hence, $\phi^*u$ is normal to $S_{\ip}$ at $x$ if $u$ is normal to $S_\Omega$ at $\xi$.  The converse statement follows since $\phi$ is a bijection.
To establish the second result, let $u$ be a 2-form which is tangent to $S_\Omega$ then $u \cdot \nabla f = 0.$
The chain rule gives us $\nabla g = (D\phi)^T (\nabla f)\circ \phi$ and by definition of the pullback, $\phi^*u = |D\phi|(D\phi)^{-1}\cdot(u\circ \phi)$, hence:
\begin{align*}
\phi^*u \cdot \nabla g &= |D\phi| (u\circ \phi)^T\cdot({D\phi^{-1}})^T\cdot(D\phi)^T\cdot[(\nabla f)\circ \phi]\\
&= |D\phi| (u^T\cdot\nabla f) \circ \phi = 0.
\end{align*}
Hence $\phi^*u$ is tangent to $S_{\ip}$.  Again, the proof of the converse follows by noting that $\phi$ is a bijection.\end{proof}

Any construction of conforming finite elements must include consideration of the traces of approximants onto inter-element boundaries. To this end, we introduce  some notation for the trace maps to the different faces of the reference pyramids.  We do not need to define traces  for the approximants in $L^2(\Omega)$.
\begin{definition}
Let $S_{i,\ip}$ be a vertical face of $\ip$. For $s=0,1,2$ define the pullback of the inclusion $S_{i,\ip}\hookrightarrow \ip$ as the trace map $\Gamma^{s}_{i,\ip}$ on  $\mathcal{U}^{s,k}(\ip)$ for all $k\in \mathbb{N}$ and $i=1..4$. We denote by $\Gamma^s_{i,\Omega}$ the corresponding  trace onto the triangular faces of  the finite pyramid $\Omega$. We similarly define the trace maps onto the base faces, that is,  $\Gamma^s_{B, \ip}$ and $\Gamma^s_{B, \Omega}$ are  the trace maps to $B_{\ip}$ and $B_\Omega$ respectively.
\end{definition}
 The consequence for us is that trace maps commute with $\phi^*$, (e.g.  $\Gamma^{s}_{i,\ip} \circ \phi^* = \phi^* \circ \Gamma^{s}_{i,\Omega}$) so results we establish on faces and edges of $\ip$ will carry over to the finite pyramid.

We can now describe the inter-element compatibility conditions to be satisfied by the traces of our approximation spaces. From \cite{monk:fe-maxwell}, we can concisely denote trace spaces on each face of the $k$th order tetrahedral and hexahedral elements by the polynomial spaces $\tau^{(s),k}$ and $\sigma^{(s),k}$ respectively. On the triangular face $S_{1,\Omega}$ and the base face $B_\Omega$, these spaces are defined as
\begin{align*}
\tau^{(0),k} (\xi,\zeta)&= P^k(\xi,\zeta)& \sigma^{(0),k}(\xi,\eta) &= Q^{k,k}(\xi,\eta) \\
\tau^{(1),k} (\xi,\zeta) &=  (P^{k-1} (\xi,\zeta))^2  \oplus S^{k,2}(\xi,\zeta)  & \sigma^{(1),k} (\xi,\eta) &= Q^{k-1,k}(\xi,\eta) \times Q^{k,k-1}(\xi,\eta)    \\
\tau^{(2),k} (\xi,\zeta) &= P^{k-1} (\xi,\zeta) & \sigma^{(2),k}(\xi,\eta)  &= Q^{k-1,k-1}(\xi,\eta) 
\end{align*}
where $S^{k,2}(\xi,\zeta) = \{w(\xi, \zeta)\in (\tilde{P}^k)^2 | (\xi - \xi_0, \zeta - \zeta_0) \cdot w = 0\}$ for some fixed $(\xi_0, \zeta_0)$.
In order to satisfy the compatibility condition {\bf (P1)}, then, we will have to enforce the constraints
 \begin{equation} \label{mega-trace-constraint}\Gamma^s_{1,\ip} u\in \tau^{(s),k}(\xi,\zeta)  \qquad \forall u\in \mathcal{U}^{(s),k}(\Omega), \forall s=0,1,2\end{equation}on the face $S_{1,\Omega}$. Analogous constraints will hold on all the other faces of the pyramid $\Omega$. 
 
The discussions above suggest the face-wise constraints which must be satisfied by any approximation spaces $\mathcal{U}^{(s),k}(\Omega)$.  However, as was demonstrated by Theorem \ref{impossible} the difficulty of  interpolation on a pyramid stems from the need to find an interpolant that match trace data on {\it all} the faces simultaneously.  This point will be discussed later.

\section{The approximation spaces $\mathcal{U}^{(s),k}(\ip)$ on the infinite pyramid}\label{infiniteapprox}
In this section we present the approximation spaces  $\mathcal{U}^{(s),k}(\ip)$ on the infinite pyramid. These will be used, via the pullback map, to construct the approximation spaces $\mathcal{U}^{(s),k}(\Omega)$ on the finite pyramid. As a preliminary step, we identify families of ``rational polynomials" on $\ip$ which will be used extensively.
We want the spaces on the finite pyramid $\Omega$ to contain all polynomials up to a specified degree.  Consider the effect of the pullback mapping $\phi$ on a polynomial of degree $k$, $p=\xix^\alpha \xiy^\beta \xiz^ \gamma \in H^1(\Omega)$, where $\alpha+\beta+\gamma = k$:  
\begin{align}\label{kweightexample}
\phi^*p = \frac{x^\alpha y^\beta z^\gamma}{(1+z)^k}.
\end{align}
From Lemma \ref{identification_theorem}, the pullback $\phi^*p\in H_w^1(\ip)$. 
This motivates our next definition: 
\begin{definition}
Let $Q^{l,m,n}(x,y,z)$ to be the space of polynomials of maximum degree $l,m,n$ in $x,y,z$ respectively.  Define the space of \emph{k-weighted tensor product polynomials} 
\begin{equation*}
Q_k^{l,m,n} (x,y,z)= \left\lbrace\frac{u}{(1+z)^k} \;:\; u \in Q^{l,m,n}(x,y,z)\right\rbrace.
\end{equation*}
It will be helpful to remember the inclusion:
\begin{align} \label{polyinc}
Q_k^{l,m,n} \subset Q_{k+1}^{l,m,n+1}.
\end{align}
Let $\PPP^n(x,y,z)$ be polynomials of maximum total degree $n$ in $(x,y,z)$ and define the space of \emph{k-weighted polynomials of degree n}
\begin{equation}\label{kweightedpdef}
\PPP_k^n (x,y,z)= \left\{\frac{u(x,y,z)}{(1+z)^k} \;:\; u(x,y,z) \in \PPP^n(x,y,z)\right\}.
\end{equation}
\end{definition}

\subsection{ $H^1_w(\ip)$-conforming approximation spaces} 
We recall from \cite{monk:fe-maxwell} that the finite element approximation space for a hexahedral element consists of polynomials of form $p = \xix^\alpha \xiy^\beta \xiz^\gamma$. From \eqref{kweightexample}, we know that $\phi^*p = \frac{x^\alpha y^\beta z^\gamma}{(1+z)^{k}} \in H_w^1(\ip),$ if $\alpha+\beta+\gamma =k.$
We might therefore expect to base an approximation space for $H_w^1(\ip)$ on the $k$-weighted space, $Q^{k,k,k}_k$.  However, there are some elements of $Q^{k,k,k}_k$ which, when pulled back to the finite pyramid, become undefined at $\xii_0=(0,0,1)$. The problem arises with elements of the form $\frac{x^ay^bz^k}{(1+z)^k}$ on the infinite pyramid. The following examples are illustrative.
\begin{example}\label{h1example}
Consider the monomial ${p_1}(x,y,z) = x$ on the infinite pyramid.  The inverse pull-back onto the finite pyramid is  $(\phi^{-1})^*{p} = \frac{\xi}{1-\zeta}$.  The limit $\lim_{\xii \rightarrow \xii_0} (\phi^{-1})^*{p}$ depends on the path by which we approach $\xii_0$.  Specifically, if we take the path $\alpha_\lambda(t) = (\lambda(1-t),0,t)$ then $\lim_{t \rightarrow 1} (\phi^{-1})^*{p}(\alpha_\lambda(t)) =\lambda$. 
\end{example}
\begin{example}
Consider the function ${p_2}(x,y,z)=\frac{z^k}{(1+z)^k}$ on the infinite pyramid. Pulled back to the finite pyramid, $(\phi^{-1})^*{p_2} =\zeta^k$. We must therefore retain ${p_2}$ in the approximation space on the infinite pyramid. 
\end{example}
\begin{lemma} \label{theoremh1}
Let $\Omega_\infty$  be the infinite pyramid described above, and $k\geq 1$ be a fixed integer.
 \begin{itemize} 
 \item Functions ${p}(x,y,z):= \frac{ x^a y^b z^{c}}{(1+z)^k} \in Q^{k,k,k-1}_k $ satisfy $p \in H^1_w(\ip)$.
\item  If $p(x,y,z) = \frac{r(x,y)z^k}{(1+z)^k}, r(x,y) \in Q^{k,k}(x,y)$, then $ \lim_{\xii \rightarrow \xii_0} (\phi^{-1})^*(p)$ is only well-defined if $r(x,y)\equiv1$.
\end{itemize}
\end{lemma}
\begin{proof}
We can verify the first statement  by using Definition 1. The second statement can be proved by contradiction, as in Example \ref{h1example}. 
\end{proof}

This result and the examples suggest the basis functions to include in a finite-dimensional approximation space for $H^1_w(\ip)$.
\begin{definition}
Let $k$ be a positive integer. We define the {\it underlying spaces $\overline{\G }(\ip)$} 
\begin{align}\label{h1spacesa}
\overline{\G }(\ip) &= Q_k^{k,k,k-1} \oplus \spn{\frac{z^k}{(1+z)^k}}.\end{align} 
\end{definition}
\begin{lemma} \label{basis_define} The rational polynomials $\left\{ \frac{x^a y^b z^c}{(1+z)^k}, 0 \leq a,b\leq k, 0\leq c\leq k-1\right\}$ and $\frac{z^k}{(1+z)^k} $ form a  basis for $\overline{\G}(\ip) $.  Moreover, $\overline{\G}(\ip) $ can be represented as 
\begin{align}\label{alth1space}
\overline{\G}(\ip) &= \{u \in Q_k^{k,k,k} : \nabla u \in Q_{k}^{k-1,k,k-1} \times Q_{k}^{k,k-1,k-1} \times Q_{k+1}^{k,k,k-1} \rbrace.
\end{align}     
\end{lemma}
\begin{proof}
The basis functions are determined by using the definition of $ \overline{\G}(\ip)$ and Lemma \ref{theoremh1}. The gradients of rational functions of the form $  \frac{x^a y^b z^c}{(1+z)^k}$ are 1-forms in $Q_{k}^{k-1,k,k-1} \times Q_{k}^{k,k-1,k-1} \times Q_{k+1}^{k,k,k-1}$. Moreover,  $\nabla \frac{z^k}{(1+z)^k} =({0},{0},{\frac{kz^{k-1}}{(1+z)^{k+1}}})^T$. The reverse inclusion  follows readily by a similar calculation. This establishes the alternative characterization of $\overline{\G}(\ip)$. 
\end{proof}

We must now constrain these spaces to obtain the approximation spaces which satisfy the compatibility constraints {\bf P1}.  This follows the discussion in Section \ref{section:rotations}, and specifically \eqref{mega-trace-constraint}. 
\begin{definition}
Let $k$ be a positive integer. We define the $k$-th order approximation spaces $\G(\ip)$:
\begin{align}\label{h1spacesb}
\G(\ip)  &= \lbrace u \in \overline{\G }(\ip)\; | \; \Gamma_{1,\ip} \in P^k_k[x,z], \text{ similarly on }S_{i,\ip},i=2,3,4\rbrace.
\end{align}
\end{definition}
%
Since we will be working in the projection-based interpolation framework while specifying internal degrees of freedom, we  define a subspace $\G _0(\ip)$, consisting of functions in $\G(\ip)$ with zero trace on the boundary of $\ip$.  Clearly,  $\G _0 (\ip)= \{x(1-x)y(1-y)zu,\; u \in Q_k^{k-2,k-2,k-2} \}$. 

%
In the Appendix, we present the shape functions in $\G(\ip)$ associated with the faces, edges and vertices of $\ip$. These are linearly independent. Moreover, the number of these functions associated with a given triangular or squareface is exactly the same as the dimension of trace spaces $\tau^{(0),k}$ or $\sigma^{(0),k}$ respectively.

\subsection{ $H_w(\curl,\ip)$-conforming  approximation spaces}

We now present the construction of the approximation space $\C(\ip)$ of $H_w(\curl,\ip)$. As before, this construction is motivated by the ultimate goal of constructing a finite element approximation space for $H_w(\curl,\Omega)$ which satisfies property ({\bf P1}). 

To satisfy the commuting diagram property we will need, at the very least, to have $\nabla \G(\ip)  \subset \C(\ip) $.  The alternate characterization of $\overline{\G }(\ip)$ in Lemma \ref{basis_define} suggests that we might consider the space $Q_{k}^{k-1,k,k-1} \times Q_{k}^{k,k-1,k-1} \times Q_{k+1}^{k,k,k-1}$ as a candidate for an approximation space for $H_w(\curl, \ip)$.   However, this space includes functions that are undefined at the point $\xii_0 = (0,0,1)$ on the finite pyramid.  We must be careful here to identify what kind of discontinuities we wish to exclude on the finite pyramid.  Firstly, we are not interested in point values of these functions, only their tangential components.  Secondly, given a particular tangent direction, $\overline{v}$ on a face of the finite pyramid, it only makes sense to consider limits taken along paths on faces tangent to $\overline{v}$. The following examples illuminate these points.
\begin{example} Consider $u = \col{y/(1+z)}{0}{0} \in Q_{k}^{k-1,k,k-1} \times Q_{k}^{k,k-1,k-1} \times Q_{k+1}^{k,k,k-1}$.  Its (inverse) pullback to the finite pyramid is, $ (\phi^{-1})^*{u} = \col{\eta/(1-\zeta)}{0}{\xi\eta/(1-\zeta)^2}$. 

 Let $\overline{v} = (0,-1,1)$ and consider the path $\alpha_\lambda(t) = (\lambda(1-t),1-t,t)$.  This path lies on the face $S_3$ for $\lambda \in [0,1]$, and $S_3$ is tangent to $\overline{v}$.  The limit of the component of $(\phi^{-1})^{*}u$ tangent to $\overline{v}$ at $\xii_0$ along the path $\alpha_\lambda$ is $\lim_{t \rightarrow 1} u(\alpha_\lambda(t))\cdot\overline{v} = \lambda$. This limit therefore depends on the path taken to approach $\xii_0$.
\end{example}
\begin{example}
Let $u = \frac{z^{k-1}}{(1+z)^{k+1}} \col{r_x z}{r_y z}{-r},\quad r \in Q^{k,k}[x,y], r_x:=\frac{\partial r}{\partial x}, r_y:=\frac{\partial r}{\partial y}$, be a 1-form defined on the infinite pyramid. 
Note that we can write $u = \nabla (\frac{r z^k}{(1+z)^{k+1}}) - \col{0}{0}{\frac{(k+1)rz^{k-1}}{(1+z)^{k+2}}}$, from which it is apparent that $ u \in H_w(\curl,\ip)$.

\end{example}
With these examples in hand, we are able to define approximation spaces for $H_w(\curl,\ip).$
\begin{definition} Let $k\geq 1$ be an integer. We define the {\it underlying space} for $H_w(\curl,\ip)$:
\begin{multline}\label{hcurlspace}
\overline{\C }(\ip) := Q_{k+1}^{k-1,k,k-1} \times Q_{k+1}^{k,k-1,k-1} \times Q_{k+1}^{k,k,k-2} \\ \oplus \left\{\frac{z^{k-1}}{(1+z)^{k+1}} \col{r_x z}{r_y z}{-r},\quad r \in Q^{k,k}[x,y]\right\}. 
\end{multline}We have again used the notation $r_x:=\frac{\partial r}{\partial x}, r_y:=\frac{\partial r}{\partial y}$.
An equivalent characterization of the underlying space $\overline{\C}(\ip)$ is given as
\begin{multline}\label{althcurlspace}
\overline{\C}(\ip) = \lbrace u \in Q_{k+1}^{k-1,k,k} \times Q_{k+1}^{k,k-1,k} \times Q_{k+1}^{k,k,k-1} : \\ \nabla \times u \in Q_{k+2}^{k,k-1,k-1} \times Q_{k+2}^{k-1,k,k-1} \times Q_{k+2}^{k-1,k-1,k} \rbrace, 
\end{multline}
\end{definition}
We now add constraints on the tangential traces,  analogous to \eqref{mega-trace-constraint}, to get the full definition of the approximation space $\C(\ip)$. Concretely,  let $n_i$ be the (outward) normal to the vertical faces $S_{i,\ip}$ of $\ip$. Then $\Gamma^1_{i,\ip} u :=  u\times n_i \vert_{S_{i,\ip}}$ for $u\in \overline{\C(\ip)}$.
\begin{definition} Let $k\geq 1$ be an integer. Define
\label{hcurlsurf}
\begin{multline}
\C(\ip)  = \bigl\lbrace u \in \overline{\C(\ip) } \; \bigg\vert \;\Gamma^1_{1,\ip} u \in (P^{k-1}_{k+1}[x,z])^2 \oplus \tilde{P}^{k-1}_{k+1}[x,1+z]\coll{1+z}{-x} \\ \text{ and similarly on }S_{i,\ip},i=2,3,4, \bigr\rbrace,
\end{multline}
where  $$\tilde{P}_{k+1}^{k-1}[x,1+z] = \frac{1}{(1+z)^{(k+1)}} \spn{x^a(1+z)^{k-1-a}, 0 \le a \le k-1 }.$$
\end{definition}
We can also identify elements in $\C(\ip)$ whose (tangential) traces vanish on $\partial \ip$. We denote the set of these as $\C_0(\ip)$. 

In the Appendix we have tabulated the edge and face shape functions for $\C(\ip)$. These are linearly independent, and are consistent along shared edges. The same will be true of the pull-backs onto the finite pyramid.

\subsection{ $H_w(\divv,\ip)$ and $L^2_w(\ip)$-conforming    approximation spaces}
Following a similar strategy to the previous sections, in this section we construct approximation spaces $\D(\ip) $ for $\HH_w(\divv, \ip)$, such that their pull-backs to the finite pyramid provide approximation spaces for $\HH(\divv,\Omega)$. Again, we want $ \curl u \in \D (\ip) ,\,\forall u\in \overline{\C}(\ip).$
Now, the curls of  functions $u \in \overline{\C }(\ip)$ satisfy  $$\nabla \times u \in  Q_{k+2}^{k,k-1,k-1} \times Q_{k+2}^{k-1,k,k-1} \times Q_{k+1}^{k-1,k-1,k-1}.$$Not all of these will have well-defined normal traces, and we must exclude these.
\begin{definition}
The underlying space for the $\HH(\divv)$-conforming element is defined as:
\begin{align}\label{hdivspace}
\begin{split}
\overline{\D }(\ip)&= Q_{k+2}^{k,k-1,k-2} \times Q_{k+2}^{k-1,k,k-2} \times Q_{k+2}^{k-1,k-1,k-1}  \\
&\quad \oplus \frac{z^{k-1}}{(1+z)^{k+2}}\col{0}{2s}{s_y (1+z)} \oplus \frac{z^{k-1}}{(1+z)^{k+2}}\col{2t}{0}{t_x (1+z)}. \\ 
\end{split}
\end{align}
Here $s(x,y) \in Q^{k-1,k}[x,y], \, s_y:=\frac{\partial s}{\partial y},$ and $ t(x,y) \in Q^{k,k-1}[x,y], t_x:=\frac{\partial t}{\partial x}$.  An alternate characterization of $\overline{\D }(\ip)$ is 
\begin{multline}\label{altdk}
\overline{\D}(\ip) = \lbrace u \in Q_{k+2}^{k,k-1,k-1} \times Q_{k+2}^{k-1,k,k-1} \times Q_{k+2}^{k-1,k-1,k} :  \nabla \cdot u \in Q_{k+3}^{k-1,k-1,k-1}\rbrace. 
\end{multline}
\end{definition}
We equip this space with  constraints on normal traces to obtain the full definition of the approximation space $\D  (\ip)$ on the infinite pyramid: 
\begin{definition} The $k$th order approximation space for $H_w(\divv,\ip)$ is
\begin{align}
\label{hdivsurf}
\D  (\ip)= \lbrace u \in \overline{\D } \; | \;\Gamma_{1,\ip}^{(2)}\in P^{k-1}_{k+2}[x,z], \text{similarly on }S_{i,\ip},i=2,3,4 \rbrace.
\end{align}
\end{definition}
Again, we can identify the 2-forms in $\D(\ip)$ with vanishing normal traces on the faces of $\ip$. We denote this set by $\D_0(\ip)$. In the Appendix, we have written down a basis for  $\D_0(\ip)$, and augmented it with shape functions for the faces.

Since we want the divergence operator to be surjective as a map from $\D (\ip)$ to the associated approximation space of $L_w^2(\ip)$, the approximation space for $L^2(\ip)$ (considered as the space of 3-forms) consists precisely of  $\divv \D(\ip) $. There is no longer any need to define an underlying space.
\begin{definition}
We define the approximation space $\Z(\ip)$ for $L^2_w(\ip)$ as
\begin{align}\label{l2space}
\Z(\ip) &= Q_{k+3}^{k-1,k-1,k-1}.
\end{align}
\end{definition}


\section{The approximation spaces  $\mathcal{U}^{(s),k}(\Omega)$ on the finite pyramid}\label{finiteapprox}
We are now readily able to define the approximation spaces for the de Rham sequence on the finite pyramid, based on the approximation spaces on the infinite pyramid $\ip$: 
\begin{definition} \label{define_pyr} Let $\Omega$ be the finite reference pyramid as defined in \eqref{pyramid_def}. Then, the $k$th order conforming subspaces on the finite pyramid $\Omega$ are
\begin{align}
 \mathcal{U}^{(s),k} (\Omega) := \set{(\phi^{-1})^*u : u \in \mathcal{U}^{(s),k}(\ip) } , \qquad s=0,1,2,3.
\end{align}
We also denote by 
\begin{align}
 \overline{\mathcal{U}^{(s),k} (\Omega)} := \set{(\phi^{-1})^*u : u \in \overline{ \mathcal{U}^{(s),k}(\ip) }} , \qquad s=0,1,2,3.
\end{align}
the underlying spaces.
\end{definition}
\begin{theorem} Let $k$ be a positive integer. The finite dimensional spaces defined in \eqref{define_pyr} satisfy:
 \begin{align}\G(\Omega) \subset H^1(\Omega), \qquad \C (\Omega)\subset H(\curl, \Omega),  \\
\D (\Omega) \subset H(\divv,\Omega), \qquad   \Z(\Omega) \subset L^2(\Omega). \end{align}
\end{theorem}
\begin{proof}
The proof follows from the definitions and properties of $\mathcal{U}^{(s),k}(\ip)$, the pull-back map $\phi$, and Lemma \ref{identification_theorem}.
\end{proof}

In the following subsections, we shall establish several useful properties of these spaces. The analysis will typically be performed for the approximation spaces on the infinite pyramid, where the basis functions are tensorial in nature, and hexahedral symmetries can be used, which allows for simple calculations in many cases. The properties of the pull-back operator will allow us to demonstrate the results on the finite pyramid.

\subsection{ $H^1(\Omega)$-conforming approximation spaces} \label{gksection}
In this subsection, we demonstrate that the grad operator is injective on $\G_0(\Omega)$, the set of bubble functions on the pyramid. \begin{lemma} \label{grad0}
Let $\G_0(\Omega)$ be the subset of $\G(\Omega)$, consisting of functions whose trace onto the faces and edges of $\Omega$ are zero. If $\nabla v=0$ for some $v\in \G_0(\Omega)$, $v\equiv 0$ on $\Omega$. 
\end{lemma}
\begin{proof} This follows from the divergence theorem.
\end{proof}

We can easily see that 
$\G_0(\Omega) = \left\{ (\phi^{-1})^*u : u \in \G_0(\ip) \right \}.$ From the remarks following \eqref{h1spacesb}, it follows that $ \dim \G_0(\Omega) = \dim \G_0(\ip) = (k-1)^3.$ Note that from the definition of $\G_0(\ip)$ and the discussion in Section \ref{section:rotations}, the face traces of functions in $\G(\Omega)$ are compatible with those of neighbouring tetrahedral and hexahedral elements. Finally, the shape functions in the Appendix show that the edge traces are well-defined, and that edge traces can be specified in  consistent manner.

\subsection{ $H(\curl,\Omega)$-conforming approximation spaces} 
We shall establish that  the $\grad$ operator maps  $\overline{\G(\Omega) }$ into $ \overline{\C}(\Omega)$. This is an important step towards showing exactness of the diagram in \ref{cd}. We then show that the curl operator is injective on a certain subspace of $\C(\Omega)$, which will be used in establishing unisolvency of the edge elements on the pyramid. We will finally demonstrate a discrete Helmholtz decomposition. Note that from the definition of $\C_0(\ip)$ and the discussion in Section \ref{section:rotations}, the face traces of functions in $\C(\Omega)$ are compatible with those of neighbouring tetrahedral and hexahedral elements.

\begin{lemma}\label{gradgkck}
The gradient operator is well defined as a map from  $ \overline{\G}(\Omega)$ into $\overline{\C}(\Omega)$.
\end{lemma}
\begin{proof}
It is easier to work on the infinite pyramid.  Recall that a basis for $\overline{\G}(\ip) $ is given by functions of the form $u_{a,b,c} = \frac{x^a y^b z^c}{(1+z)^k}$, where $a,b$ and $c$ are integers and $a \in [0,k]$, $b \in [0,k]$ and $c \in[0,k-1]$ or $u_{a,b,c}= \frac{z^k}{(1+z)^k}$.  We will show that the gradients of each of these functions lie in $\overline{\C}(\Omega)$.
The result is trivial for $c=0$.  For $c \ge 1$,
\begin{align*}
\nabla u_{a,b,c} = \frac{1}{(1+z)^{k+1}}\col{a x^{a-1}y^{b} (z^{c+1}+z^c)}{b x^{a}y^{b-1} (z^{c+1}+z^c)}{x^{a}y^{b}((c-k)z^c + c z^{c-1})}.
\end{align*}
 If $c \le k-2$ then $\nabla u_{a,b,c} \in Q_{k+1}^{k-1,k,k-1} \times Q_{k+1}^{k,k-1,k-1} \times Q_{k+1}^{k,k,k-2}$.  In the case $c = k-1$, we can let $r = x^a y^b$ in \eqref{hcurlspace} and then the remainder 
\begin{align}
\nabla u_{a,b,c} - \frac{z^{k-1}}{(1+z)^{k+1}}\col{r_x z}{r_y z}{-r} = \frac{1}{(1+z)^{k+1}}\col{a x^{a-1}y^{b}z^{k-1}}{b x^{a}y^{b-1} z^{k-1}}{c x^{a}y^{b}z^{k-2}},
\end{align}
which is in $Q_{k+1}^{k-1,k,k-1} \times Q_{k+1}^{k,k-1,k-1} \times Q_{k+1}^{k,k,k-2}$.  Finally, if $c=k$ then choosing $r=1$  in \eqref{hcurlspace} suffices. 
Now use the definition of $\overline{\mathcal{U}^{(s),k}}(\Omega)$ in terms of the  inverse  pull-back of functions in $\overline{\mathcal{U}^{(s),k}}(\ip)$, and the commutativity of the grad with the pull-backs, to conclude the result.  
\end{proof}

Note that the previous result also follows immediately from the (unproven) equivalent characterisations of the underlying spaces, \eqref{alth1space} and \eqref{althcurlspace}.  
An important subset of $\C(\Omega)$ is  the functions  with vanishing tangential traces.
\begin{definition} \label{defck0curl}
Define $\C_0(\Omega)$ to be  the subspace of functions in $\C(\Omega)$ with zero tangential component on the boundary of $\Omega$.  
\end{definition}

From Lemma \ref{tracemaps}, we know that if $u\in \C(\Omega)$ has zero tangential traces on a particular face or edge of $\Omega$, then its pullback to $\ip$ will have zero tangential traces on the associated face or edge. This allows us to characterize $\C_0(\Omega)$.\begin{lemma}\label{c0omega}
Functions in $\C_0(\Omega)$ can be represented as $(\phi^{-1})^*(u)$, where $ u \in \C_0(\ip)$ have the form 
\begin{align}\label{c0omegadef}
u=  \col{y(1-y)zq_1}{x(1-x)zq_2}{x(1-x)y(1-y)q_3} + \frac{z^{k-1}}{(1+z)^{k+1}} \col{r_x z}{r_y z}{-r},
\end{align}
where $q \in Q_{k+1}^{k-1,k-2,k-2} \times Q_{k+1}^{k-2,k-1,k-2} \times Q_{k+1}^{k-2,k-2,k-2}$ and $r = x(1-x)y(1-y)\rho,\; \rho \in Q^{k-2,k-2}[x,y]$.  We have denoted $r_x:=\frac{\partial r}{\partial x}, r_y:=\frac{\partial r}{\partial y}$.
\end{lemma}
\begin{proof} It is easily verified that the functions $u$ in \eqref{c0omegadef} have zero tangential traces on the edges and faces of $\ip$, and therefore their inverse pullbacks $(\phi^{-1})^*(u)$ belong to $\C_0(\Omega)$. Note also that
$$  \dim \C_0(\Omega)  = \dim \C_0(\ip) =  k(k-1)^2+k(k-1)^2+(k-1)^3+(k-1)^2 = 3k(k-1)^2.$$
\end{proof}

The curl operator has a non-empty null space in $\C_0(\Omega)$, consisting of gradients. We can precisely characterize the complement of the gradients in $\C_0(\Omega)$.
\begin{definition}
Define  $\C_{0,\curl}(\Omega) \subset \C_0(\Omega)$ as
$$\C_{0,\curl}(\Omega):= \left\{ v \vert  v = (\phi^{-1})^*u, u \in \C_{0,\curl} (\ip)\right\}$$,  where $\C_{0,\curl}(\ip) \subset \C_0(\ip)$ consists of functions $u$ of the form
\begin{align} \label{ck0curl}
u=\col{y(1-y)z q_1}{x(1-x)z q_2}{x(1-x)y(1-y)\rho},\end{align}
with$ q_1 \in Q_{k+1}^{k-1,k-2,k-2},  q_2 \in Q_{k+1}^{k-2,k-1,k-2},\rho \in Q_{k+1}^{k-2,k-2}[x,y].$
\end{definition}
We now show that $\C_{0,\curl}(\Omega)$ contains no gradients.
\begin{lemma}\label{curlinj}
Let $\C_{0,\curl}(\Omega)$ be defined as above. 
Then $\C_{0,\curl}(\Omega)\subset \C_0(\Omega)$, and the curl operator is injective on $\C_{0,\curl}(\Omega)$. In other words, 
$ \grad \G_0(\Omega)\cap \C_{0,\curl}(\Omega) = \{ 0\}$.
\end{lemma}
\begin{proof} The set inclusion $\C_{0,\curl}(\Omega)\subset \C_0(\Omega)$ follows by the definitions of $\C_{0,curl}(\Omega)$ and $\C_0(\Omega)$. To see that the curl operator is injective on $\C_{0,\curl}(\Omega)$, we first show that the curl operator is injective on $\C_{0,\curl}(\ip)$. The argument proceeds by contradiction.

If $k=1$ then $\C_{0,\curl}(\ip)$ is empty.  Assume $k \ge 2$ and let $u \in \C_{0,\curl}(\ip)$ be as in \eqref{ck0curl}.  Let either $\rho$ or $q_2$ not equal to zero and write $\rho = \frac{r(x,y)}{(1+z)^{k+1}}$, $r \in Q^{k-2,k-2}(x,y)$.  Suppose that $\nabla \times u =0$.  From the $x$-component, we obtain
\begin{align*}
\frac{1}{(1+z)^{k+1}}\pd{}{y}\left(y(1-y)r\right) - \pd{}{z}\left(zq_2\right) = 0.
\end{align*}
There is no $z$-dependence in $r$ so we can factorise $q_2 = f(z)g(x,y)$, where $f \in P^{k-2}(z)$ satisfies
\begin{align*}
\frac{d}{dz}\frac{zf(z)}{(1+z)^{k+1}} = \frac{1}{(1+z)^{k+1}}.
\end{align*}
This is impossible, and so $\rho= q_2 = 0$.  A similar consideration of the $y$-component shows that $q_1=0$.
We have just established that the curl operator is injective on $\C_{0,\curl}(\ip)$. Since the pullback and curl commute, the curl is injective on $\C_{0,\curl}(\Omega)$.
\end{proof}

We can now state a discrete Helmholtz decomposition for $\C_0(\Omega)$: 
\begin{lemma} The discrete approximation space $\C_0(\Omega) \subset \HH(\curl, \Omega)$ of functions with vanishing tangential traces on $\partial \Omega$ admits a Helmholtz decomposition. That is, if $v \in \C_0(\Omega)$, we can write $v = \nabla q + w$ with $q \in \G_0(\Omega)$ and $w \in \C_{0,\curl}(\Omega)$.  \label{curldecomp}
\end{lemma}
\begin{proof}
If $q\in \G_0(\Omega)$, it has zero trace on all the faces and edges of $\Omega$. Therefore, the tangential components of $\nabla q$ are also zero on the faces and edges. We already know that $\grad$ maps $\overline{\G}(\Omega)$ into $\overline{\C}(\Omega)$ from 
Lemma \ref{gradgkck}, and so it is clear that $\grad$ maps $\G_0(\Omega)$ into $\C_0(\Omega)$. Injectivity of this map follows from Lemma \ref{grad0}
Now we count dimensions. From Section \ref{gksection} we saw that $\dim \G_0(\Omega) = (k-1)^3$, and from Lemma \ref{curlinj},  
$$ \dim \C_{0,\curl}(\Omega) = \dim \C_{0,\curl}(\ip) = k(k-1)^2+k(k-1)^2+(k-1)^2 = (2k+1)(k-1)^2.$$
From the same lemma, we know $\grad \G_0 (\Omega)\cap \C_{0,\curl} (\Omega)= 0$.  Both of these are subspaces of $\C_0(\Omega)$. So, 
$$\dim \left\{ \grad \G_0 (\Omega) \cup \C_{0,\curl} (\Omega)\right\}  =  (2k+1)(k-1)^2 + (k-1)^3 = 3k(k-1)^2$$
which is the dimension of $\C_0(\Omega)$.  Hence $\C_0(\Omega) =  \grad \G_0 (\Omega) \oplus \C_{0,\curl} (\Omega).$ 
\end{proof}

\subsection{ $H(\divv,\Omega)$-conforming approximation spaces } 

In this subsection we shall establish that $\curl \overline{\C}(\Omega) \subset \overline{\D}(\Omega)$. We then show that the div operator is injective on a certain subspace of $\D(\Omega)$. We finally demonstrate a decomposition of this discrete space.
\begin{lemma}\label{curlckdk}
The curl operator maps elements of $\overline{\C}(\Omega)$ into $\overline{\D} (\Omega)$.  
\end{lemma}
The proof of this lemma is a calculation similar to the one in Lemma \ref{gradgkck}, and is omitted here.

We now need to identify elements of $\D(\Omega)$ which have vanishing normal traces on the faces of the finite pyramid. Denote these by $\D_0(\Omega)$. 
%
From Lemma \ref{tracemaps}, we know that if $\Gamma^2_{i,\Omega} u =0$ for some $u\in \D(\Omega)$, then the pull-back $\Gamma^2_{i,\ip} \phi^*u =0$ on the associated face of $\ip$. This allows us to characterize $\D_0(\Omega)$ easily.
\begin{lemma} Functions in $\D_0(\Omega)$ can be represented as $(\phi^{-1})^*(u)$, where $ u \in \D_0(\ip)$ have the form 
\begin{align}\label{dk0}
\begin{split}
&\frac{z^{k-1}}{(1+z)^{k+2}}\col{2t}{2s}{(1+z)\left(s_y+t_x\right)} +\col{x(1-x)\chi_1}{y(1-y)\chi_2}{z\chi_3},\end{split}
\end{align}
where $ s=y(1-y)\sigma,\; t=x(1-x)\tau,$
with $\chi_1 \in Q_{k+2}^{k-2,k-1,k-2},\; \chi_2 \in Q_{k+2}^{k-1,k-2,k-2},$ $ \chi_3 \in Q_{k+2}^{k-1,k-1,k-2}$, $\sigma \in Q^{k-1,k-2}(x,y), s_y:=\frac{\partial s}{\partial y}$, and $ \tau \in Q^{k-2,k-1}(x,y), t_x:=\frac{\partial t}{\partial x}$.
\end{lemma}
\begin{proof} It is easily verified that functions of the form \eqref{dk0} have vanishing normal components on the faces $S_{i,\ip}$ of the infinite pyramid; their (inverse) pullbacks to the finite pyramid will thus have vanishing normal components on the faces $S_{i,\Omega}$ of $\Omega$. 
\end{proof}

We note also that
\begin{align*} \dim \D_0(\Omega) &= \dim \D_0(\ip) \\ &= k(k-1)^2 + k(k-1)^2 + k^2(k-1) + k(k-1) + k(k-1) \\
&=3k^3-3k^2. \hspace{2in} \end{align*}

We now present a subspace of $\D_0(\Omega)$ on which the divergence operator will be injective.
\begin{definition} Define $\D_{0,\divv} (\Omega) := \{ v \vert v = (\phi^{-1})^*(u), u \in \D_{0,\divv}(\ip)\}$ where
\begin{equation}\label{divg}  \D_{0,\divv}(\ip):=\spnn \{ \frac{z^{k-1}}{(1+z)^{k+2}}\col{r_y+2t}{r_x + 2s}{(1+z)(r_{xy}+s_y+t_x)} \}\oplus \spnn\{ \col{0}{0}{ z\chi_3} \} \end{equation}
and where 
$ r(x,y) = x(1-x) y(1-y)p(x,y), p\in Q^{k-2,k-2},  t = x(1-x)\tilde t, \tilde{t} \in P^{k-2}(x)$, $s=y(1-y)\tilde{s}, \tilde{s} \in P^{k-2}(y), $ and $\chi_3 \in  Q_{k+2}^{k-1,k-1,k-2}.$ Again, the subscripts denote partial differentiation.
\end{definition}
\begin{lemma} \label{divinj}
The divergence operator is injective on  $\D_{0,\divv}(\Omega)$.
\end{lemma}
\begin{proof} We shall first show that the divergence operator is injective on $\D_{0,\divv}(\ip)$. 
Let $u$ be as in \eqref{divg}. If $\nabla \cdot u=0$, then 
\[ 0=\nabla \cdot u = \frac{(k-1)z^{k-2}}{(1+z)^{k+2}}(r_{xy}+t_x+s_y) + \frac{\partial}{\partial z} (z \chi_3).\]
We  factorize $\chi_3 = \sum_{i=0}^{k-2} \frac{z^i}{(1+z)^{k+2}} q_i(x,y)$ and  compare coefficients of like powers of $z$. Since $r,t$ and $s$ have no dependence on $z$, we obtain
\begin{align*} 0&= \frac{(k-1)z^{k-2}(r_{xy}+t_x+s_y)}{(1+z)^{k+2} } +  \frac{d}{dz}\left( \sum_{i=0}^{k-2}\frac{z^{i+1}q_i(x,y)}{(1+z)^{k+2}} \right) \\
&=\frac{(k-1)z^{k-2}(r_{xy}+t_x+s_y)}{(1+z)^{k+2} } + \left( \sum_{i=0}^{k-2}\frac{z^{i+1}(i-k-1) + (1+i)z^i}{(1+z)^{k+3}}q_i(x,y) \right) .\end{align*}
This is impossible,  unless
\[ (r_{xy}+t_x+s_y) =0, \quad q_i(x,y)=0.\] 
However, $t$ only depends on $x$, and $s$ only depends on $y$. From the form of $r$, it must be that $r\equiv 0 \equiv t \equiv s$. Therefore, $\nabla u \not= 0 $ for any non-zero $u\in \D_{0,\divv}(\ip)$. Using the properties of the pullback operator, $\nabla\cdot v =0 \Rightarrow v\equiv 0$ for all $v \in \D_{0,\divv}(\Omega)$. 
The desired result on $\Omega$ will follow by the properties of the pullback operator $\phi$ and the commutativity of $\phi$ with the divergence.
 \end{proof}
 It is easy to see that \[ \dim \D_{0,\divv}(\Omega) = \dim \D_{0,\divv}(\ip) = (k-1)^2+2(k-1)+ k^2(k-1) = k^3-1. \]
Just as in the previous section, we can use Lemma (\ref{divinj}) to exhibit a convenient decomposition of the discrete approximation space.
\begin{lemma} 
\label{divdecomp}
Any $v \in \D_0(\Omega)$ can be decomposed as $v = \nabla \times w_1+ w_2$ with $w_1\in \C_{0,\curl}(\Omega)$, $w_2 \in \D_{0,\divv}(\Omega)$.
\end{lemma}
\begin{proof}
Lemma \ref{curlckdk} tells us that the curl operator maps $\overline{\C}(\Omega)$ into $\overline{\D}(\Omega)$.  Observe that if the tangential components of $v$ are zero on some surface then the component of $\nabla \times v$ that is normal to the surface will also be zero and so the curl operator maps $\C_{0,\curl}(\Omega)$ into $\D_0(\Omega)$.  By Lemma \ref{curlinj} we know that this mapping is injective. 

By construction,  $\D_{0,\divv}(\Omega)$ is  a subset of $\D_0(\Omega)$ and by lemma \ref{divinj}, $\nabla \cdot w \ne 0$ for all $w \in \D_{0,\divv}(\Omega)$.  Hence $\D_{0,\divv}(\Omega) \cap \C_{0,\curl}(\Omega)$ is empty.
We now count dimensions. We established in the proof of Lemma (\ref{divinj}) that  $\D_{0,\divv}(\Omega)$ has dimension $k^3-1$ and from the previous section we know $\C_{0,\curl}(\Omega)$ has dimension $2k^3 - 3k^2 +1$.  Thus, 
\[ \dim (\curl \C_{0,\curl}(\Omega) \cup \D_{0,\divv}(\Omega) ) = 3k^3-3k^2 = \dim  \D_0(\Omega),\] which shows that
$ \D_0(\Omega) = \curl \C_{0,\curl}(\Omega) \oplus \D_{0,\divv}(\Omega) $.
This establishes the desired decomposition.
\end{proof}

\subsection{ $L^2(\Omega)$-conforming approximation spaces} We note that the dimension of $\Z(\Omega) = \dim \Z(\ip) = \dim(Q^{k-1,k-1,k-1}_{k+3}) = k^3$.  
It is a straightforward matter to determine that the divergence operator is well defined as a map from $\overline{\D}(\Omega)$ to $\Z(\Omega)$. We record the result here in a lemma.
\begin{lemma}\label{divdkzk}
The divergence operator maps elements of $\overline{\D}(\Omega)$ into  $\Z(\Omega)$. 
\end{lemma}
\begin{lemma} Any element $u\in \Z(\Omega)$ can be written uniquely as 
$$ u = \nabla \cdot w+ \lambda, \qquad w \in \D_{0,\divv}(\Omega),\; \lambda \in \mathbb{R}.$$\label{l2decomp}
\end{lemma}
\begin{proof} From Lemma \ref{divdkzk}, we know that $\divv  \D_{0,\divv}(\Omega) \subset \Z(\Omega)$. We also know that the constants belong to $\Z(\Omega)$. Now, $\dim (\divv  \D_{0,\divv}(\Omega) ) = k^3-1$, which is one less than the dimension of $\Z(\Omega)$. 
Now, suppose we could find $w \in \divv  \D_{0,\divv}(\Omega)$ so that $\nabla w =1$ on $\Omega$. By definition of $\D_{0,\divv}(\Omega)$, we know that $w$ has zero normal components on the faces of $\Omega$. From the divergence theorem, this is impossible. Hence, we have shown that the constants are not contained in $\divv  \D_{0,\divv}(\Omega)$, and therefore 
$ \divv  \D_{0,\divv}(\Omega) \oplus {\mathbb{R}} = \Z(\Omega).$ This completes the proof.
\end{proof}

We finish this subsection with an important component of the proof that our elements satisfy property {\bf P1}. 
 \begin{lemma}\label{tracesokay}
The spaces of traces of the approximation spaces, $\mathcal{U}^{(s),k}(\Omega)$ on the faces of the pyramid are the same as those of the corresponding tetrahedral and hexahedral elements.  Specifically  $\Gamma^s_{i,\Omega}\mathcal{U}^{(s),k}(\Omega) = \tau^{(s),k}$ and $\Gamma^s_{B,\Omega} \mathcal{U}^{(s),k}(\Omega)= \sigma^{(s),k}$.
\end{lemma}
\begin{proof}
In the Appendix, we collect shape functions in Tables  \ref{tab:basisfunctions0}, \ref{tab:basisfunctions1} and \ref{tab:basisfunctions2}  of the approximation spaces $\mathcal{U}^{(s),k}(\Omega)$ for $s=0,1$ and 2 respectively.  It can also be easily (though tediously) verified that the traces of these shape functions on each face span the corresponding trace space from the tetrahedral and hexahedral elements.  This demonstrates that $\Gamma^s_{i,\Omega}\mathcal{U}^{(s),k}(\Omega) \supseteq \tau^{(s),k}(S_{i,\Omega})$ for $i=1,2,3,4$ and $\Gamma^s_{B,\Omega}\mathcal{U}^{(s),k}(\Omega) \supseteq \sigma^{(s),k}(B)$.  

Set equality is seen by examining the infinite pyramid case. By construction, if $u\in \mathcal{U}^{s,k}(\ip)$, then its trace $\Gamma^{(s)}_{1,\ip} u $ on the vertical face $S_{1,\ip}$   lies in  $P_k^k$, $(P^{k-1}_{k+1}[x,z])^2 \oplus \tilde{P}^{k-1}_{k+1}[x,1+z]\coll{1+z}{-x}$ or $P^{k-1}_{k+2}$ for $s = 0,1,2$ respectively. This means that $dim\, (\Gamma^{(s)}_{i,\Omega} \mathcal{U}^{(s),k}(\Omega) \leq dim\, \tau^{(s),k}(S_{i,\Omega)}$, which establishes that $(\Gamma^{(s)}_{i,\Omega} \mathcal{U}^{(s),k}(\Omega) = \tau^{(s),k}(S_{i,\Omega)}$.  Also, elements of  $\Gamma^{(s)}_{i,\Omega}\mathcal{U}^{s,k}(\Omega)$ consist of the pullbacks of functions in $\Gamma^s_{i,\ip}\mathcal{U}^{(s),k}(\ip)$. Therefore, 
 rotational symmetry means similar statements hold for the other faces as well.  Finally, the dimension of $\Gamma^s_{B,\Omega}\mathcal{U}^{(s),k}(\Omega)$ is equal to that of $\sigma^{(s),k}(B)$ and so $\Gamma^s_{B,\Omega}\mathcal{U}^{(s),k}(\Omega) = \sigma^{(s),k}(B)$

\end{proof}
The implication of Lemma \ref{tracesokay} is important: the spaces $\mathcal{U}^{(s),k}(\Omega)$ allow for full compatibility  of relevant traces with well-known tetrahedral and hexahedral finite elements, across interelement boundaries. This should allow for the seamless integration of pyramidal elements into a hybrid mesh consisting of tetrahedra and hexahedra.

\section{First order elements on the pyramid}\label{firstorder}
The approximation spaces $\mathcal{U}^{(s),k}(\Omega)$ constructed above include the elements presented by \cc{Gradinaru and Hiptmair }{}\cite{hiptmair:pyramid} as the special case $k=1$.  To demonstrate this, we will map the basis functions presented in that paper onto the infinite pyramid, and demonstrate that these (pulled-back) elements belong to $\mathcal{U}^{(s),k}(\ip)$. The properties of  the pullback then allow us to conclude the set inclusions on the finite pyramid. The reason for  this indirect approach is the tensorial nature of the approximation spaces on $\ip$, which makes it easier to examine basis functions.

\begin{itemize}
\item {\bf The lowest-order $H^1(\Omega)$ element}: 
The basis functions for the $H^1(\Omega)$ element given in \cc{equation 3.2 of \cite{hiptmair:pyramid}}{\cite[equation 3.2]{hiptmair:pyramid}} are denoted $\pi_i,\; i=1..5$.  Set $\tilde{\pi}_i = \phi^*\pi_i$. 
\begin{equation*}\label{g1basis}
\begin{split}
\tilde{\pi}_1 &= \frac{(\xx-1)(\xyy-1)}{1+\xz}, \quad
\tilde{\pi}_2 = \frac{\xx(\xyy-1)}{1+\xz}, \quad
\tilde{\pi}_3 = \frac{(\xx-1)\xyy}{1+\xz}, \quad
\tilde{\pi}_4 = \frac{\xx\xyy}{1+\xz},\\
\tilde{\pi}_5 &= \frac{\xz}{1+\xz}. \quad
\end{split}
\end{equation*}
 It is clear that $\tilde{\pi}_i \in \mathcal{U}^{(0), 1}(\ip).$
\item {\bf The lowest-order $\HH(\curl,\Omega)$ element:}
We proceed as in the $H^1(\Omega)$ case. Set $\tilde{\gamma}_i = \phi^*\gamma_i$ where the $\gamma_i,\; i=1..8$ are the basis functions for the curl-conforming element in \cite{hiptmair:pyramid}\footnote{There are minor typographical errors in \cite{hiptmair:pyramid} for two of the one-forms.  Based on the preceding calculations in that paper, the correct expressions are
\begin{equation*}
\gamma_6 = \col{-z + \frac{yz}{1-z}}{\frac{xz}{1-z}}{x-\frac{xy}{1-z}+\frac{xyz}{(1-z)^2}}, \quad \gamma_7 = \col{\frac{yz}{1-z}}{-z +\frac{xz}{1-z}}{y-\frac{xy}{1-z}+\frac{xyz}{(1-z)^2}}
\end{equation*}

}:
\begin{equation}\label{c1basis}
\begin{split}
\tilde{\gamma}_1 &= \frac{1}{(1+\xz)^2}\col{1- \xyy}{0}{0}, \quad 
\tilde{\gamma}_2 =  \frac{1}{(1+\xz)^2}\col{0}{\xx}{0},\quad 
\tilde{\gamma}_3 = \frac{1}{(1+\xz)^2}\col{\xyy}{0}{0}, \\
\tilde{\gamma}_4& =  \frac{1}{(1+\xz)^2}\col{0}{1- \xx}{0},\quad 
\tilde{\gamma}_5 =  \frac{1}{(1+\xz)^2}\col{\xz(1-\xyy)}{\xz(1-\xx)}{(1-\xyy)(1-\xx)}, \quad 
\tilde{\gamma}_6 = \frac{1}{(1+\xz)^2}\col{\xz(\xyy-1)}{\xz\xx}{\xx(1-\xyy)}, \quad \\ 
\tilde{\gamma}_7 &= \frac{1}{(1+\xz)^2}\col{\xz\xyy}{\xz(\xx-1)}{\xyy(1-\xx)}, \quad 
\tilde{\gamma}_8 = \frac{1}{(1+\xz)^2}\col{-\xz\xyy}{-\xz\xx}{\xx\xyy}
\end{split}
\end{equation}
These are also the pullbacks of the basis functions for the first order curl conforming element given by \cc{Graglia et al. }{}\cite{graglia:highorderpyramid}. Note that these are all edge shape functions.  It is easy to see that $\tilde{\gamma}_i $ are  shape functions specified in the previous section for $H_w(\curl, \ip)$ with $k=1$.

\item {\bf The lowest-order $\HH(\divv,\Omega)$ element:}
Set $\tilde{\zeta}_i = \phi^*\zeta_i$, where $\zeta_i,\; i=1..5$ are the divergence-conforming basis functions 
\begin{equation}\label{d1basis}
\begin{split}
\tilde{\zeta}_1 &= \frac{1}{(1+\xz)^3}\col{0}{2(\xyy-1)}{\xz}, \quad
\tilde{\zeta}_2 = \frac{1}{(1+\xz)^3}\col{2(\xx-1)}{0}{\xz}, \\
\tilde{\zeta}_3 &= \frac{1}{(1+\xz)^3}\col{2\xx}{0}{\xz}, \quad
\tilde{\zeta}_4 = \frac{1}{(1+\xz)^3}\col{0}{2\xyy}{\xz}, 
\tilde{\zeta}_5 = \frac{1}{(1+\xz)^3}\col{0}{0}{-1}.
\end{split}
\end{equation}
\end{itemize}

For completeness, we note that $\mathcal{U}^{(3),1}(\Omega)$ consists of the constants, which map to multiples of $\frac{1}{(1+\xz)^4}$ on the infinite pyramid. The above collections of functions are consistent with the definitions \eqref{h1spacesa}, \eqref{hcurlspace}, \eqref{hdivspace} and \eqref{l2space}.

%% file: highorderpyramid-part2.tex


\section{The degrees of freedom $\Sigma^{(s),k}$ and unisolvency}\label{section:unisolvency}

We now define degrees of freedom $\Sigma^{(s),k}$ which are linearly independent and unisolvent for the finite element approximation spaces $\mathcal{U}^{(s),k}(\Omega)$. Our construction is based on the premise of ``patching'' as discussed in \cite{hiptmair:pyramid}:  ``the traces of discrete differential forms onto any interelement boundary (a (n-1)-face) have to be unique and they have to be fixed by the degrees of freedom associated with that face". This means the exterior degrees of freedom for $\Omega$ must be identical with those of neighbouring tetrahedra or hexahedra. 
Thus, to satisfy property {\bf P1}, we insist that the degrees of freedom are the same on interelement boundaries (vertices, edges and faces) as those from neighboring tetrahedra and hexahedra.  Another important consideration is locality.  \cc{The authors in \cite{hiptmair:pyramid}}{\cite{hiptmair:pyramid}} correctly identify that: ``expressions for integrals on edges contained on a face $S_{i,\Omega}$ should only depend on the degrees of freedom on that face''; addressing this challenge reveals the difficulty of treating a pyramid as a degenerate finite hexahedral element.  In our case, the degrees are chosen to be local \emph{ab initio}, but the challenge is to prove unisolvency.

In this section we use the same exterior degrees of freedom as specified in \cite{monk:fe-maxwell}. We show that these are indeed dual to the exterior shape functions specified in the Appendix. We then have to specify degrees of freedom for the remaining objects in the approximation space; for these we use the projection-based degrees of freedom as in \cc{Demkowicz }{}\cite{demkowicz:derham}. We finally show that the set of degrees of freedom are unisolvent. 
Throughout this and the subsequent sections, if $P$ is some finite-dimensional vector space, we will use the notation  $\cB{P}$ to denote an arbitrary basis. 

\subsection{ $H^1$-conforming element} 
In order to fully describe the $H^1-$conforming finite element on a pyramid, we need to specify 4 classes of functionals which form a dual set to the approximating basis functions: vertex, edge, face and volume degrees of freedom.  We call the set of these functions $\Sigma^{(0),k}$, and then show that $(\Omega, \G(\Omega), \Sigma^{(0),k})$ is a conforming and unisolvent element for $H^1(\Omega)$. We shall follow the presentation in \cc{Chapter 5 of \cite{monk:fe-maxwell}}{\cite[chapter 5]{monk:fe-maxwell}}.   

Depending on $k$ not all of the degrees of freedom will be needed. We explicitly design the vertex, edge and face classes of these degrees of freedom to match those of tetrahedral or hexahedral elements.  In order that the function evaluations be well-defined, let $p\in H^{3/2+\epsilon}(\Omega)$.  
\begin{enumerate}
\item {\it Vertex degrees of freedom:} let $v_i$, $ i=1..5$ be the vertices of the finite pyramid. Then $M_V $ is the set of vertex degrees of freedom $m_{v_i}$ where
$$ m_{v_i}(p) := p(v_i),  i=1..5.$$
These are identical to the vertex degrees of freedom on tetrahedral or hexahedral elements.
\item {\it Edge degrees of freedom:} these are given by the set $M_{E}$ of functionals of the form
\begin{subequations}
\begin{equation} m_{e,q}(p) := \int_{e} pq \,ds, \quad q \in \cB{P^{k-2}(e)},\qquad \text{for each edge, $e$}. \label{h1edge}\end{equation}

There are  $k-1$ linearly independent functionals $m_{e,q}$ for each of the eight edges $e \in E$. The form of these degrees of freedom is the same for ``vertical'' edges, $e_i$, and base edges, $b_i$.
Again, these are identical to the edge degrees of freedom on tetrahedral or hexahedral elements. If $k<2$ these degrees of freedom are not used.
\item {\it Face degrees of freedom: }the degrees of freedom on the triangular faces, $M_{S}$ correspond to those on the faces of tetrahedral elements.  They have the form:
\begin{equation} m_{S_i,q}(p) = \int_{\Gamma_i} pq \,dA, \quad q \in \cB{P^{k-3}(S_{i,\Omega})}, \quad i=1..4.\label{h1face}\end{equation} 
There are $(k-1)(k-2)/2$ such degrees for each face.

The degrees of freedom on the base face, $M_B$ correspond to those for hexahedral elements:
\begin{equation}m_{B,l}(p) = \int_B pq \,dA,  \quad q \in \cB{Q^{k-2,k-2}(B)}. \end{equation}
There are $(k-1)^2$ such degrees.  The face degrees of freedom are $M_F = M_S \cup M_B$. If $k<2$ these degrees of freedom are not used. 

\item {\it Volume degrees of freedom:} denote by $\G_0(\Omega)$ the subset of $\G(\Omega)$ with zero boundary traces. Then the volume degrees of freedom are given by 
\begin{equation}M_\Omega := \left\{p \mapsto \int_{\Omega} \nabla p\cdot \nabla q \,dV, \quad q \in \cB{\G_0(\Omega)} \right\}\label{h1degrees}.\end{equation} 
The dimension of $\G_0(\Omega)$ is $(k-1)^3$.  If $k<2$ these degrees of freedom are not used.\end{subequations}
\end{enumerate}
The set of all degrees of freedom for $s=0$ is $\Sigma^{(0),k}:= M_V \cup M_E\cup M_F\cup M_\Omega$.
We can now state the major conformance and unisolvency result:
\begin{theorem} \label{h1conformingunisolvent} The element $(\Omega, \G(\Omega), \Sigma^{(0), k})$ is $H^1$-conforming and unisolvent.\end{theorem}
\begin{proof} 
To show that this element is conforming, we need to establish that  the vertex, edge and face degrees of freedom of $p\in \G(\Omega)$ vanish on a face of the pyramid, if and only if $p\equiv 0$ on that face.  By Lemma \ref{tracesokay} the trace $\Gamma^{0}_{i,\Omega}p$ to the triangular face $S_{i,\Omega}$ lies in  $\tau^{(0),k}$. The trace $\Gamma^{0}_{B,\Omega}p$ lies in $\sigma^{(0),k}$. Now, we have chosen the degrees of freedom so that on each each face they are also identical to to those of the corresponding (conforming) tetrahedral or hexahedral element. The vanishing of the external degrees of freedom associated with a face therefore implies that $p\equiv 0$ on that face, see \cc{Lemmas 5.47 and 6.9 of \cite{monk:fe-maxwell}}{\cite[lemmas 5.47 and 6.9]{monk:fe-maxwell}}

For unisolvency we need to show that for any vector $(u_i) \in \RR^{\dim{\Sigma^{(0),k}}}$, there exists a unique element $u \in \G(\Omega)$ with $m_i(u) = u_i \;\forall m_i \in \Sigma^{(0),k}$.  We first observe that $\dim{\Sigma^{(0),k}} = k^3+3k+1= \dim  \G(\Omega)$ and so uniqueness implies existence, i.e. we need to show that if {\it all} the degrees of freedom of $p\in \G(\Omega)$ vanish, then $p\equiv 0$ on $\Omega$.  We have just seen that the vanishing of the external degrees of freedom implies $p=0$ on $\partial \Omega$ and hence $p\in \G_0(\Omega)$. The vanishing of the volume degrees of freedom implies that
$$ \int_\Omega \nabla p\cdot \nabla q \,dV =0 , \quad \forall q \in \G_0(\Omega).$$ Hence $\int_\Omega |\nabla p|^2 \,dV=0,$ from which we easily see that $p\equiv 0$.
\end{proof}

\subsection{$H(\curl)$-conforming element} 
A curl-conforming pyramidal element is defined by the triple $(\Omega, \C(\Omega), \Sigma^{(1), k})$ where the degrees of freedom $\Sigma^{(1), k}$ are associated with the edges, faces,  and volume of the pyramid. Again, we follow the presentation of \cite{monk:fe-maxwell}: let $ t$ be a unit tangent vector along the edge $e$, $\nu$ be the normal to a given face, and let $u\in H^r(\curl, \Omega)$ be smooth enough so that the following functionals are well-defined:

\begin{subequations}
\begin{enumerate}
\item {\it Edge degrees of freedom}:
\begin{equation}
M_E := \set{ u \mapsto \int_{e} u \cdot t q \,ds, \quad \forall q \in \cB{P^{k-1}(e_i)}  \quad \forall e \in E }
\end{equation} 
\item {\it Face degrees of freedom}: here we must differentiate between the triangular and square faces of the pyramid. 
On the triangular faces, we specify face degrees of freedom which are identical to those for tetrahedral elements:
\begin{align}
M_{S} := \set{ u \mapsto \int_{S_{i, \Omega}} u \cdot  q \,dA, \quad \forall q \in \cB{T} \quad i = 1..4}
\end{align}
where $T = \set{q \in (P^{k-2}(S_{i,\Omega})^3 \; | \; q \cdot \nu = 0}$ and on $B$, the degrees of freedom are identical to those for hexahedral elements:
\begin{align}
M_B := \set{ u \mapsto \int_B u\cdot q \,dA, \quad \forall q \in \cB{Q^{k-2,k-1}(B) \times Q^{k-1,k-2}(B) }}. 
\end{align}
The class of face degrees of freedom is then $M_F = M_{S} \cup M_B$.
\item {\it Volume degrees of freedom}:
here we must specify the degrees of freedom associated with ``gradient bubbles'' $\nabla \G_0(\Omega)$ and ``curl bubbles'' $\C_{0,\curl}(\Omega)$. 
\begin{align}
M_{\Omega}^{\grad} &:= \set{ u \mapsto  \int_{\Omega} u \cdot \nabla q \,dV, \quad \forall q \in \cB{\G_0(\Omega)} }, \label{hcurlvol1}\\
M_{\Omega}^{\curl} &:= \set{ u \mapsto \int_{\Omega} \nabla \times u \cdot \nabla \times v \,dV, \quad \forall v \in \cB{\C_{0,\curl}(\Omega)} } \label{hcurlvol2}
\end{align}
The volume degrees are then $M_\Omega:= M_{\Omega}^{\grad}\cup M_{\Omega}^{\curl}$ 
\end{enumerate}
\end{subequations}
We must demonstrate that the finite element $(\Omega, \C(\Omega), \Sigma^{(1), k})$ is indeed curl-conforming and that specifying the degrees of freedom for a $u \in \C(\Omega)$ uniquely specifies $u$. This is the content of the next theorem:
\begin{theorem}\label{curlunisolvence}
The element $(\Omega, \C(\Omega), \Sigma^{(1), k})$ is curl-conforming and unisolvent.\end{theorem}
\begin{proof}
By an analogous argument to that given for the  $s=0$ case in Theorem \ref{h1conformingunisolvent} we see that the vanishing of the external degrees for any $u \in \C(\Omega)$ implies that $u \in \C_0(\Omega)$ and thus that the element is conforming.  All that remains is to show that if  $u\in \C_0(\Omega)$ and all the volume degrees also vanish, then $u \equiv 0$. From Lemma \ref{curldecomp} we can write
$$ u =  \nabla q' + v', \qquad q' \in \G_0(\Omega), v' \in \C_{0,\curl}(\Omega).$$ 
Since the gradient-bubble degrees of freedom, $M_{\Omega}^{\grad}(u)$ vanish, we have
$$ \int_{\Omega} (v' + \nabla q')  \cdot \nabla q \,dV =0 , \quad \forall q \in \G_0(\Omega), \Rightarrow \int_\Omega |\nabla q'|^2 =0.$$
This allows us to conclude $\nabla q'=0$. Moreover, since the curl-bubble degrees of freedom $M_{\Omega}^{\curl}(u)$ also vanish, we have
$$\int_{\Omega} \nabla \times (v'+\nabla q') \cdot \nabla \times v \,dV =0, \quad \forall v \in \C_{0,\curl}(\Omega) \Rightarrow \int_{\Omega} |\nabla \times v'|^2 \,dV=0.$$ Since the curl map was injective on $\C_{0,\curl}(\Omega)$, we know that $v'=0$. This establishes unisolvency.
\end{proof}

\begin{figure}
\centering
\includegraphics{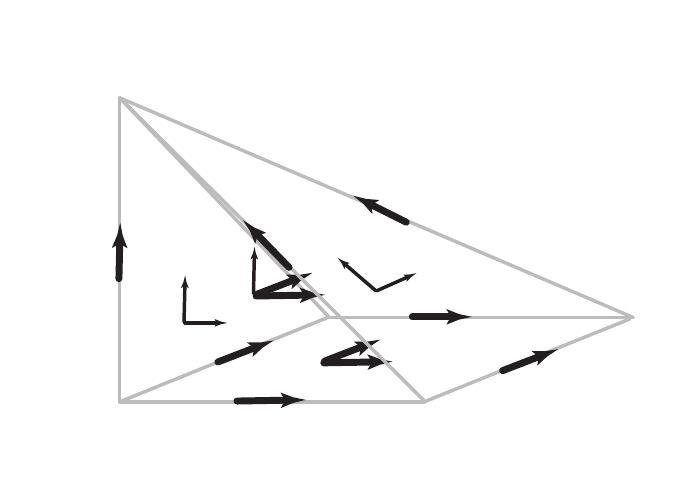}
\caption{A representation of the curl degrees of freedom for $k=2$.  The degrees solely associated with the two rear triangular faces have been omitted.  Bold arrows indicate two degrees of freedom.  $\mathcal{U}^{(0),2}_0$ contributes one volume degree and $\C_{0,\curl}$ contributes four (two pairs).}
\label{fig:curldegrees}
\end{figure}

\subsection{$H(\divv)$-conforming element} 
By now the strategy of defining a conforming element using the space $\mathcal{U}^{(s),k}$ is familiar: we define exterior degrees of freedom to ensure conformancy, and use a Helmholtz-like decomposition of the approximation space to ensure unisolvency. For the triple $(\Omega, \D(\Omega), \Sigma^{(2), k})$, we define the degrees of freedom by specifying the face and volume degrees:
\begin{subequations} 
\begin{enumerate}
\item {\it Face degrees of freedom}: we have to specify separate degrees of freedom on the triangular and square faces. On the triangular faces $S_{i,\Omega}, i=1..4$, we specify degrees of freedom $M_S$ in terms of the basis functions $q$ of $ (P^{k-1}(S_{i,\Omega})$. On the base face $B$, we specify the face degrees of freedom $M_B$ in terms of the basis function $q$ of  $Q^{k-1,k-1}(B)$.
\begin{align}
M_S &:= \set{u \mapsto  \int_{S_{i,\Omega}} u \cdot  \nu q\, dA, \quad \forall q \in \cB{(P^{k-1}(S_{i,\Omega})},\quad i=1..4}\\
M_B &:= \set{u \mapsto \int_B u\cdot \nu q\, dA, \quad \forall q \in \cB{Q^{k-1,k-1}(B)}}.
\end{align} 
The set of face degrees of freedom are then $M_F = M_S \cup M_B$.
\item {\it Volume degrees of freedom}: $M_\Omega:= M_\Omega^{\curl} \cup M_\Omega^{\divv}$ where
\begin{align}
M_\Omega^{\curl} &:= \set{u \mapsto \int_{\Omega} u \cdot \nabla \times v\, dV,\quad \forall v \in \cB{\C_{0,\curl} (\Omega)}} \label{hdivvol1},\\
M_\Omega^{\divv} &:= \set{u \mapsto \int_{\Omega} \nabla \cdot u\, \nabla \cdot v\, dV,\quad \forall v \in \cB{\D_{0,\divv} (\Omega)}} \label{hdivdegrees}.
\end{align}
\end{enumerate}
Again, $\Sigma^{(2),k}:= M_F\cup M_\Omega$.
\end{subequations} 
\begin{theorem} The finite element triple 
 $(\Omega, \D(\Omega), \Sigma^{(2), k})$ is divergence-conforming and unisolvent.\end{theorem}
\begin{proof}
Conformance follows by an argument similar to that for Theorems \ref{h1conformingunisolvent} and \ref{curlunisolvence}. For unisolvency, if all the degrees of freedom for a given $u\in \D(\Omega)$ vanish, then we must show that $u\equiv 0$. Now, since the element is conforming, we know that vanishing face degrees of freedom means $u\in \D_0(\Omega)$. 

From Lemma \ref{divdecomp}, $u\in \D_0(\Omega)$ can be written as $u =  \nabla \times w_1+ w_2$ with $w_1\in \C_{0,\curl}(\Omega)$, $w_2 \in \D_{0,\divv}(\Omega)$. The vanishing of the $M_\Omega^{\curl}(u)$ and $M_\Omega^{\divv}(u)$ degrees of freedom implies that $  \nabla \times w_1=0, \divv w_2=0$. Now, the curl operator is injective on  $\C_{0,\curl} (\Omega)$ from Lemma \ref{curlinj}, and so $w_1=0$. The div operator is injective on $\D_{0,\divv}(\Omega)$ from Lemma \ref{divinj}, and so $w_2=0$. This establishes the result.\end{proof}

\subsection{ $L^2$-conforming element} 
Functions in $L^2(\Omega)$ do not have well-defined traces on $\partial \Omega$, so we only need to specify volume degrees of freedom to completely define the finite element triple $(\Omega, \Z(\Omega), \Sigma^{(3),k}).$ The volume degrees specify the contribution from the ``divergence bubble'' and the constants
\begin{subequations}
\begin{align}
&M_\Omega := \set{p \mapsto \int_{\Omega} p \nabla \cdot v \,dV,\quad \forall v \in \cB{\D_{0,\divv}(\Omega)}}, \label{l2omegadegs}\\
&M_1(p) =  \set{p \mapsto \int_{\Omega} p \,dV }. \label{l2constdegs}
\end{align} This specifies $\Sigma^{(3),k}:= M_\Omega \cup M_1$. 
Unisolvency follows immediately by using Lemma \ref{l2decomp}.
\end{subequations}


\section{Interpolation and exact sequence property}
We have now constructed approximation subspaces $\mathcal{U}^{(s),k}$ for $H^1(\Omega), H(\curl,\Omega), H(\divv,\Omega)$ and $L^2(\Omega)$. During the process of construction, we saw that $d \overline{\mathcal{U}}^{(s),k} (\Omega) \subset \overline{\mathcal{U}}^{(s+1),k} (\Omega)$ for $s=0,1,2$. In this section, we define interpolation operators $\Pi^{(s)}$ so that the finite elements satisfy a commuting diagram property. This will enable us to  show that in fact the approximation space $\mathcal{U}^{(s),k}(\Omega)$ satisfy an exact sequence property. The degrees of freedom induce an interpolation operator on each element.  We have to be  careful about choosing the spaces that we are able to interpolate; for example, the vertex degrees for the $H^1$-conforming element require us to take point values, which are not defined for a general $H^1(\Omega)$ function.  Details of the regularity required for the external degrees can be found in \cite{monk:fe-maxwell}.  The problem is discussed for projection-based interpolation in \cite{demkowicz:projection}.  For our purposes it is enough to know that it is possible to choose $r>1$ such that all the degrees of freedom are well defined on the spaces $H^r(\Omega)$, $\HH^{r-1}(\curl, \Omega)$,  and $\HH^{r-1}(\divv, \Omega)$.  
The sets of degrees of freedom then induce interpolation operators in the expected way. 
\begin{definition} Let $k\in \mathbb{N}$ be given, and let $u$ be an $s-$form, $s=0,1,2,3$ which possesses enough regularity such that the degrees of freedom $\Sigma^{(s),k}(u)$ are well-defined. We define the local interpolation operator $\Pi^{(s)}$ by requiring that  $\Pi^{(s)}(u) \in \mathcal{U}^{(s),k}(\Omega)$ and for all degrees of freedom $m\in \Sigma^{(s),k}$, 
\begin{align}\label{interpdef}
 m(u) = m(\Pi^{(s)}u).
 \end{align}
\end{definition}
The interpolation operator  is well-defined, since the $\Sigma^{(s),k}$ are unisolvent. It is is local on  the faces, edges and vertices of $\Omega$, and agrees with the choice for high-degree elements presented in \cite{monk:fe-maxwell}.  Therefore, the construction of a global interpolant on a mesh which includes pyramidal elements will be simple.  The volume degrees of freedom are reminiscent of, and inspired by, the projection-based interpolation framework of \cite{demkowicz:projection}. Providing optimal $hp$ estimates of the interpolation error in this framework is technical, and relies on the use of a basis-preserving extension operator. We leave this for future work.

Equipped with these interpolation operators, the finite elements satisfy a commuting diagram property: 
\begin{theorem}Let $r>0$ be chosen so that the interpolation operators $\Pi^{(s)}$ are well-defined. Then the diagram 
\begin{equation} \label{cd2}
\begin{CD}
H^r(\Omega) @>d>> \HH^{r-1}(\curl, \Omega) @>d>> \HH^{r-1}(\divv, \Omega) @>d>> H^{r-1}(\Omega) \\
@V \Pi^{(0)} VV @V \Pi^{(1)} VV @V \Pi^{(2)} VV @V \Pi^{(3)} VV \\
\G (\Omega) @>d>> \C (\Omega) @>d>> \D (\Omega) @>d>> \Z(\Omega)
\end{CD}
\end{equation}
commutes.
\end{theorem}

\begin{proof} For each $s=0,1,2$, we have to show that $d\Pi^{s}p = \Pi^{s+1}dp$ for any $s-$form, $p$.  This is equivalent to showing that $\Pi^{(s+1)} d(p-\Pi^{(s)}p) = 0$, which, in turn is equivalent to showing that
\begin{align}\label{degreescommute}
 \forall m \in \Sigma^{(s+1),k}, \quad m(d(p-\Pi^{(s)}p) =0,.
 \end{align} 

We split the proof by considering the exterior degrees of freedom seperately.  For each $s=0,1$, the external degrees of freedom are identical to those stated in \cite{monk:fe-maxwell}. Therefore we can adopt components of the proofs of commutativity from  \cite{nedelec:mixed1980,monk:fe-maxwell} to see that the 
$ m(d(p-\Pi^{(s)}p) =0$ for each exterior degree of freedom, $m\in \Sigma^{(s+1),k}$, $s=0,1$. There are no external degrees of freedom in $\Sigma^{(3),k}$. 

We still need to demonstrate \eqref{degreescommute} for the volume degrees of freedom in $\Sigma^{(s+1),k}$. The argument follows that of \cite{demkowicz:projection}. Let $s=0$. There are two classes of volume degrees of freedom in $\Sigma^{(1),k}$.  The first is given in \eqref{hcurlvol1}.  Let $m_v \in M_\Omega^{\curl}$ be a degree of freedom associated with the test function $v \in \C_{0,\curl}(\Omega)$
\begin{align*}
m_v(d (p - \Pi^{(0)} p)) &= \int_{\Omega} \nabla \times \nabla ( p - \Pi^{(0)}\nabla p) \cdot \nabla \times v\, dV = 0.
\end{align*}
 The second type of volume degree is given in \eqref{hcurlvol2}.  Let $m_q \in M_\Omega^{\grad}$ be the degree of freedom associated with some $q \in \G_0(\Omega)$. Then
\begin{align} \label{nontrivialdegrees}
m_q(d (p -\Pi^{(0)} p)) &= \int_{\Omega} \nabla (p -\Pi^{(0)}p) \cdot \nabla q\, dV = 0 
\end{align}
because of the definition of the interpolation operator, \eqref{interpdef} and the $H^1$ volume degrees of freedom, \eqref{h1degrees}.   Here the important point is that the same function spaces is used  in each of these sets of degrees of freedom.  The proof for $s=1$  follows from a similar argument, this time using the equivalence of \eqref{hdivvol1} and \eqref{hcurlvol2} to deal with the homogenous divergence-free part.  

For $s=2$, the degrees, $M_\Omega$ given in \eqref{l2omegadegs} can be dealt with in the same fashion as \eqref{nontrivialdegrees}.  For the final degree of freeom, $M_1$, given in \eqref{l2constdegs}, we note that 
\begin{align*}
\int_\Omega \nabla \cdot (p-\Pi^{(2)}p) = \int_{\partial \Omega} (p-\Pi^{(2)}p) \cdot \nu dS = 0
\end{align*}
because we have already established the commutativity of the external degrees and the test functions used for the external degrees include constants on each face.
%
%
\end{proof}

\begin{theorem} 
The following sequence is exact
\begin{equation}
\begin{CD}
\RR @>>> \G(\Omega) @>\nabla>> \C(\Omega) @>\nabla\times>> \D(\Omega)@>\nabla\cdot>> \Z(\Omega)@>>> 0.
\end{CD}
\end{equation}
\end{theorem}
\begin{proof} We need to show the inclusions  $d\,\mathcal{U}^{(s),k}(\Omega) \subset \mathcal{U}^{(s+1),k} (\Omega)$ for $s=0,1,2$ and the property if  $u$ is an $s+1$ form with $du=0$, then $u=dv$ for some $v\in\mathcal{U}^{(s),k}(\Omega)$.

By the definitions, \eqref{h1spacesa}, \eqref{hcurlspace}, \eqref{hdivspace} and \eqref{l2space}, we see that $d\overline{\mathcal{U}^{(s),k}}(\Omega) \subset \overline{\mathcal{U}^{(s+1),k}} (\Omega)$ for $s=0,1,2$. By Theorem \ref{tracesokay} it follows that the face restrictions inherit the exact sequence property for tetrahedral and hexahedral elements, so that $d\,\mathcal{U}^{(s),k}(\Omega) \subset \mathcal{U}^{(s+1),k} (\Omega)$.
%
%

To show the second property, which is equivalent to demonstrating the existence of discrete potentials, we shall use Theorem \ref{cd}. 
First let  $s=0$, and suppose $u\in \mathcal{U}^{(1),k}(\Omega)$ satisfies $\nabla \times u=0$. Then there is a continuous $v\in H^1(\Omega)$ such that $u= \nabla v$. Using the commuting diagram property, $ u= \Pi^{(1)} u = \Pi^{(1)} \nabla v = \nabla \Pi^{(0)}v $, and thus $u$ is derivable from a discrete potential.  The argument for $s=1$ and $s=2$ is identical.
\end{proof}


\section{Polynomial approximation property}
We now need to show that our approximation spaces $\mathcal{U}^{(s),k}(\Omega)$ allow for high-degree approximation. Concretely, given any desired degree $q\in \mathbb{N}$, we need to demonstrate that we can choose $k$ so that polynomials of degree $q$ are contained in $\mathcal{U}^{(s),k}(\Omega)$.  We start with the $L^2$-conforming element.

\begin{lemma} \label{zpolyprop}
The $L^2$-conforming element exactly interpolates all polynomials up to degree $k-1$.  That is, $P^{k-1} (\Omega)\subset \Z(\Omega)$.
\end{lemma}
\begin{proof}
A basis for $P^{k-1}(\Omega)$ is given by functions of the form \[u=\xi^a \eta^b (1-\zeta)^c\] where $a, b, c$ are non-negative integers and $a + b + c \le k-1$.  Using the pullback formula, we see that
\begin{equation*}
\phi^*u=\frac{x^a y^b (1+z)^{k-1-(c+a+b)}}{(1+z)^{k+3} },
\end{equation*}
which is in $Q_{k+3}^{k-1,k-1,k-1} = \Z(\Omega)$.
\end{proof}

\begin{lemma}
The $H(\divv)$-conforming element exactly interpolates all polynomials up to degree $k-1$.  That is, $\PPP^{k-1} \subset \D(\Omega)$.
\end{lemma}
\begin{proof}
A basis for $\PPP^{k-1}$ is given by functions of the form:
\small
\begin{align}\label{ppk1basis}
\col{\xi^{a_1} \eta^{b_1} (1-\zeta)^{c_1}}{0}{0},\quad \col{0}{\xi^{a_2} \eta^{b_2} (1-\zeta)^{c_2}}{0},\quad \col{0}{0}{\xi^{a_3} \eta^{b_3}. (1-\zeta)^{c_3}}
\end{align}
\normalsize
where the $a_i, b_i, c_i$ are non-negative integers and $a_i + b_i + c_i \le k-1$.  Pullback these functions to the infinite pyramid to get:
\small
\begin{equation} \label{hdivpolys}
\col{\dfrac{x^{a_1} y^{b_1} (1+z)^{\overline{c_1}}}{(1+z)^{k+2} } }{0}{0},\;\col{0}{\dfrac{x^{a_2} y^{b_2} (1+z)^{\overline{c_2}}}{(1+z)^{k+2} } }{0},\;\frac{x^{a_3} y^{b_3} (1+z)^{\overline{c_3}}}{(1+z)^{k+2} } \col{x}{y}{1+z}.
\end{equation}
\normalsize
Here we have written $\overline{c_i} = k-1-(c_i+a_i+b_i)$. The constraint $a_i + b_i + c_i \le k-1$ ensures that if $\uu$ is as in \eqref{ppk1basis}, then $\phi^*\uu \in Q_{k+2}^{k,k-1,k-1} \times Q_{k+2}^{k-1,k,k-1} \times Q_{k+2}^{k-1,k-1,k}$.  Moreover, since divergence commutes with pullback, so $\nabla \cdot \phi^*\uu = \phi^* \nabla \cdot \uu$. Now $\uu \in \PPP^{k-1} \Rightarrow \nabla \cdot \uu \in \PPP^{k-2}$ and so by Lemma~\ref{zpolyprop}, $\phi^*\nabla \cdot \uu \in Q_{k+3}^{k-1,k-1,k-1}$. We have thus established that $\phi^*(u) \in  \overline{\D}(\ip)$, where we used the  characterization of the underlying space \eqref{hcurlspace}. 

Now, since $u$ is a polynomial 2-form, its normal trace onto a triangular face of $\Omega$  will be a polynomial of the same or less degree,  and hence the surface constraints in \eqref{hdivsurf} will be satisfied automatically.  Hence $\phi^*\uu \in \D(\ip)$, which means that $u \in \D(\Omega)$.
\end{proof}

The existence of polynomials in the $\HH(\curl)$-conforming element may be proved in a similar manner:  
\begin{lemma} \label{cpolyprop}
The $\HH(\curl)$-conforming element exactly interpolates all polynomials up to degree $k-1$. That is, $\PPP^{k-1} \subset \C(\Omega)$
\end{lemma}
\begin{proof}
Take basis functions for $\PPP^{k-1}$ as in \eqref{ppk1basis}.  The pullbacks of these 1-forms to the infinite pyramid are:
\small
\begin{align}
\dfrac{x^{a_1} y^{b_1} (1+z)^{\overline{c_1}}}{(1+z)^{k+1}}\col{z+1}{0}{-x},\; \dfrac{x^{a_2} y^{b_2} (1+z)^{\overline{c_2}}}{(1+z)^{k+1} }\col{0}{z+1}{-y},\; \col{0}{0}{\dfrac{x^{a_3} y^{b_3} (1+z)^{\overline{c_3}}}{(1+z)^{k+1} } } . \label{hcurlpolys}
\end{align}
\normalsize
The constraint on the $a_i, b_i, c_i$ ensures that these functions are all members of $Q_{k+1}^{k-1,k,k} \times Q_{k+1}^{k,k-1,k} \times Q_{k+1}^{k,k,k-1}$. We then use the commutativity of the curl operator with the pull-back,  the previous lemma, and the fact that the tangential traces $\Gamma^{(1)}_{i,\Omega}u$ for polynomial 1-forms $u$ satisfy the surface constraints of \eqref{hcurlsurf},  shows that the functions in \eqref{hcurlpolys} are in fact in $\C(\ip).$
\end{proof}

For the $H^1$-conforming element, we gain an extra degree in the polynomials (in fact, there are some polynomials of degree $k$ present in $\C(\Omega)$ and $\D(\Omega)$, but not all of them).
\begin{lemma} \label{gpolyprop}
The $H^1$-conforming element exactly interpolates all polynomials up to degree $k$. That is, $P^{k} \subset \G(\Omega)$
\end{lemma}
\begin{proof}
Let $p = \xi^a \eta^b \zeta^c$, $a,b,c$ be non-negative integers and $a + b + c \le k$.  
\begin{align}
\phi^*p = \frac{x^a y^b z^c (1+z)^{k-(a+b+c)}}{(1+z)^k},
\end{align}
If $a+b \not=0$ it is clear that $\phi^*p \in Q_k^{k,k,k-1}$. On the other hand, if $a+b = 0$, we obtain $\phi^*p \in \{ \frac{z^k}{(1+z)^k}\}$. Therefore, polynomial zero-forms $p$ of the form \eqref{gpolyprop} satisfy $\phi^*p \in Q_k^{k,k,k-1} \oplus \frac{z^k}{(1+z)^k} = \overline{\G}(\Omega)$, as required.  
Arguments similar to the previous cases demonstrate the inclusion $\phi^*p\in \G(\ip)$, and hence $p\in \G(\Omega)$.
\end{proof}

\section{Conclusion}
We have shown that the finite element approximation spaces $\mathcal{U}^s_k(\Omega)$ equipped with the external degrees of freedom from \cite{monk:fe-maxwell} and projection-based interpolation for the internal degrees of freedom are unisolvent and satisfy a commuting diagram property.  All the $k$th order spaces include the complete family of polynomials of degree $k-1$ and the $H^1$-conforming space includes all the degree $k$ polynomials too.  

These finite element spaces are based on rational basis functions. It is not surprising that arguments which rely on the polynomial or highly differentiable nature of regular finite element spaces will fail in the current situation. In upcoming work we present a careful analysis of quadrature errors for these approximation spaces.

%% file: Combined-tables-appendix.tex
\section{Appendix: Shape functions}
In Tables \ref{tab:basisfunctions0}, \ref{tab:basisfunctions1} and \ref{tab:basisfunctions2}, we present shape functions for $\mathcal{U}^{(s),k}(\Omega)$ for each $s=0,1,2$. This is not a hierarchical construction, and no attention has been paid to the conditioning of any resulting stiffness matrices.  

\begin{table}[htbp]
   \centering
   \begin{tabular}{@{} >{\small}l >{\small}l >{\small}p{6cm}@{}} 
   \multicolumn{3}{c}{Representative shape functions for 0-forms on a pyramid. }\\   \hline
         
      Infinite Pyramid    & Finite Pyramid& Comments \\ \hline

          $\dfrac{(1-x)(1-y)}{(1+z)^k}$ & $\dfrac{\alpha_\xi\alpha_\eta}{(1-\zeta)^{2-k}}$ & Vertex  function associated with vertex $v_1$.\\
          $\dfrac{z^k}{(1+z)^k}$& $\zeta^k$ &  Vertex function associated with vertex $v_5$.\\
              $\dfrac{(1-x)(1-y)z^a}{(1+z)^k}$  &$\dfrac{\alpha_\xi\alpha_\eta\zeta^a}{(1-\zeta)^{2+a-k}}$    &  Edge functions associated with edge $e_1$,  $1\leq a\leq k-1$. \\
$\dfrac{(1-y)(1-x)x^a}{(1+z)^k}$& $\dfrac{\alpha_\xi\alpha_\eta\xi^a}{(1-\zeta)^{2+a - k}}$   &   Edge  functions associated with base edge $b_1$, $1\leq a\leq k-1$. \\
      $ \dfrac{(1-x)(1-y)x^az^b}{(1+z)^k}$    & $\dfrac{\alpha_\xi\alpha_\eta\xi^a\zeta^b}{(1-\zeta)^{2+a+b-k}}$  &  Face shape functions associated with triangular face $S_{1,\Omega}$,  $1\leq a, b, a+b\leq k-1$.\\
      $ \dfrac{(1-x)(1-y)x^ay^b}{(1+z)^k}$    & $\dfrac{\alpha_\xi\alpha_\eta\xi^a\eta^b}{(1-\zeta)^{2+a+b-k}}$  &  Face shape functions associated with base face $B$, $1\leq a, b \leq k-1$.\\
     
   $  \dfrac{x(1-x)y(1-y)z x^ay^bz^c}{(1+z)^{k}} $&$\dfrac{ \xi^{a+1}\eta^{b+1}\zeta^{c+1}\alpha_\xi\alpha_\eta}{(1-\zeta)^{5+a+b+c - k}}$&   Volume shape functions, $0\leq a,b,c\leq k-2.$ \\
 \end{tabular}
   \caption{Shape functions on a pyramid. Since the approximation space $\mathcal{U}^{(0),k}(\ip)$ is invariant under the rotation, $R_\infty:\ip\rightarrow \ip$,  it is only necessary to demonstrate shape functions for a representative base vertex, vertical edge, base edge and vertical face. Then, using \eqref{rotation} and the subsequent remarks, the inverse pullback of these to the finite pyramid will also be invariant under the rotation $R$. Note that $\alpha_\xi:= (1-\zeta-\xi)$ and $\alpha_\eta:=(1-\zeta-\eta)$}
   \label{tab:basisfunctions0}
  \end{table}

\begin{table}[htbp]
   \centering
   \begin{tabular}{@{} >{\small}l>{\small}l>{\small}m{4cm}@{}} 

      \multicolumn{3}{>{\small}c}{Representative shape functions for 1-forms on a pyramid. }\\   \hline

       Infinite Pyramid    & Finite Pyramid& Comments \\ \hline
 $ \col{0}{0}{\frac{(x-1)(y-1)(1+z)^{c}}{(1+z)^{k+1}   }}$&$\col{0}{0}{\frac{\alpha_\eta\alpha_\xi}{(1-\zeta)^{3+c-k}}}$& Edge functions associated with $e_1$,  $0\leq c\leq k-1$.\\
 
 $\frac{1}{(1+z)^{k+1}} \col{x^c(1-y)}{0}{0}$& $ \frac{\xi^c\alpha_\eta}{(1-\zeta)^{1+c-k}}\col{1}{0}{\frac{\xi}{1-\zeta}}$ &Edge functions associated with base edge $b_1$, $0\leq c\leq k-1$.\\ 

$ \frac{z(1-y)}{(1+z)^{k+1}}\col{x^{a}(1+z)^{c}}{0}{0}$ &$\frac{\zeta\alpha_\eta\xi^a}{(1-\zeta)^{2+c+a-k}}\col{1}{0}{\frac{\xi}{1-\zeta}}$&Face functions associated with triangular face $S_{1,\Omega}$, $ 0 \leq a,c, \quad c+a \le k-2$. \\
 
$\frac{x(1-x)(1-y)}{(1+z)^{k+1}}\col{0}{0}{x^{a}z^{c}} $& $ \col{0}{0}{\frac{\xi^{1+a}\alpha_\xi\alpha_\eta\zeta^c}{(1-\zeta)^{4+a+c-k}  }}$&Face functions for $S_{1,\Omega}$, $0 \leq c,a, \quad c+a \leq k-3 $.\\

$\frac{(1-y)(1-x)x^{a}(1+z)^{k-a-2}}{(1+z)^{k+1}}\col{z}{0}{-x} $& $ \col{\alpha_\eta \alpha_\xi\xi^a \zeta}{0}{ -\alpha_\eta\alpha_\xi\xi^{a+1}}$&Face functions for $S_{1,\Omega}$, $0 \le a \le k-2 $.\\

$\frac{1}{(1+z)^{k+1}}\col{y(1-y)x^ay^b}{0}{0}$ &$ \frac{\alpha_\eta\xi^a\eta^{b+1}}{(1-\zeta)^{2+a+b-k}}\col{\eta}{0}{\frac{\xi}{1-\zeta}}$& Face shape functions for base face $B_{\Omega}$,  $a \le k-1, b_1 \le k-2 $.\\

$\frac{1}{(1+z)^{k+1}}\col{0}{x(1-x)x^{a}y^{b}}{0}$ &$\frac{\alpha_\zeta\xi^{a+1}\eta^b}{(1-\zeta)^{2+a+b-k}}\col{0}{\xi}{\frac{\eta}{1-\zeta}}$& Face shape functions for $B_{\Omega}$, $a \le k-2, b\le k-1$.\\

$\frac{y \left( 1-y \right) z\,x^a y^b z^c}{(1+z)^{k+1}}\col{1}{0}{0}$&
$\frac{\alpha_\eta  \xi^a \eta^{b+1}\zeta^{1+c}}{(1-\zeta)^{3+a+b+c-k}}\col{1 }{0}{\frac{\xi}{1-\zeta} }$
& Volume shape functions, $0\leq a\leq k-1$, $0\leq b,c\leq k-2$.\\

$\frac{x \left( 1-x \right) z x^a y^b z^c\,}{(1+z)^{k+1}}\col{0}{1}{0}$&
$\frac{\xi^{a+1}\eta^b\zeta^{1+c}\alpha_\xi}{(1-\zeta)^{3+a+b+c-k}}\col{0 }{ 1}{ \frac{\eta}{1-\zeta} }$
& Volume shape functions, $0\leq b\leq k-1$, $0\leq a,c\leq k-2$.\\

$\frac{x \left( 1-x \right) y \left( 1-y \right) x^a y^b z^c}{(1+z)^{k+1}}\col{0}{0}{1}$&
$\frac{\xi^{1+a}\eta^{1+b}\zeta^c\alpha_\xi\alpha_\eta\tilde{q}}{(1-\zeta)^{3+a+b+c-k}}\col{0 }{0 }{ \frac{1}{(1-\zeta)^{2}}}$
& Volume shape functions, $0\leq a,b,c \leq k-2$.\\

$\displaystyle \,\frac{z^{k-1}}{(1+z)^{k+1}} \col{{\frac{\partial r}{\partial x}}\,z}{\frac{\partial r}{\partial y} \,z}{-r} $&
$\displaystyle \zeta^{k-1} \col{\zeta\frac{\partial \tilde{r}}{\partial \xi}(1-\zeta) }{\zeta \frac{\partial \tilde{r}}{\partial \eta}(1-\zeta) }{- \tilde{r} +\zeta(\xi\frac{\partial \tilde{r}}{\partial \xi}+ \eta \frac{\partial \tilde{r}}{\partial \eta} )} $& \multirow{2}{4cm}{Volume shape functions, $0 \leq a,b\leq k-2$.}\\
$\quad r = x(1-x)y(1-y)x^ay^b,\;$&$\quad \tilde{r} = \frac{\xi^{1+a}\alpha_\xi\eta^{1+b}\alpha_\eta}{(1-\zeta)^{a+b+4}}$ \\\\ \hline

   \end{tabular}
   \caption{Curl-conforming shape functions  on a pyramid. Since the approximation space $\mathcal{U}^{(1),k}(\ip)$ is invariant under the rotation, $R_\infty:\ip\rightarrow \ip$,  it is only necessary to demonstrate shape functions for a representative base vertex, vertical edge, base edge and vertical face. Then, using  \eqref{rotation}  and the subsequent remarks, the inverse pullback of these to the finite pyramid will also be invariant under the rotation $R$. There are three distinct types of shape functions for the vertical faces, two for the base face, and four for the volume. Note that $\alpha_\xi:= (1-\zeta-\xi), \alpha_\eta:=(1-\zeta-\eta)$}
   \label{tab:basisfunctions1}
  \end{table}

\begin{table}[htbp]
   \centering
   \begin{tabular}{@{} >{\small}l>{\small}l>{\small}m{4cm}@{}} 

 \multicolumn{3}{>{\small}c}{Representative shape functions for 2-forms on a pyramid. }\\   \hline          
      Infinite Pyramid    & Finite Pyramid& Comments \\ \hline
$ \frac{1}{(1+z)^{k+2}} \col{0}{2(1-y)x^az^b}{-z^k} $ &$\zeta^k\col{\frac{\xi}{1-\zeta}}{\frac{\eta}{(1-\zeta)}}{-1} + \col{0}{\frac{2\alpha_\eta\xi^a\eta^b}{(1-\zeta)^{2+a+b-k}  }}{0} $ & Face shape functions associated with $S_{1,\Omega}$, $ a,b \ge 0,\; a+b \le k-1$.\\

$\frac{1}{(1+z)^{k+2}} \col{0}{0}{x^ay^b}$ & $\displaystyle  \left( 1-\zeta  \right)  ^{k-a-b-1} \col{{\xi}^{1+a}{\eta}^{b} }{{\xi}^{a}{\eta}^{b+1}}{ {\xi}^{a}{\eta}^{b}} $& Base face shape functions, $  0 \le a,b \le k-1$.\\


$\frac{z^{k-1}}{(1+z)^{k+2}}\col{2t}{0}{(1+z)(t_x)}$&$ \zeta^{k-1}\col{2\tilde{t}}{0}{0} + \zeta^{k-1}\tilde{t}_x\col{-\frac{\xi}{1-\zeta}}{\frac{-\eta}{1-\zeta}}{1}$&   \multirow{2}{4cm}{Volume shape functions, $0\leq a\leq k-2, 0\leq b\leq k-1$.}\\
$\quad t=x(1-x)x^ay^b$&$\quad \tilde{t}=\xi^{1+a}\alpha_\xi\eta^b(1-\zeta)^{-a-b-2}$&\\

$\frac{z^{k-1}}{(1+z)^{k+2}}\col{0}{2s}{(1+z)(s_y)}$&$ \zeta^{k-1}\col{0}{2\tilde{s}}{0} + \zeta^{k-1}\tilde{s}_y\col{-\frac{\xi}{1-\zeta}}{\frac{-\eta}{1-\zeta}}{1}$&   \multirow{2}{4cm}{Volume shape functions, $0\leq a\leq k-1, 0\leq b\leq k-2$.}  \\
$\quad s=y(1-y)x^ay^b$&$\quad \tilde{s}=\xi^{a}\alpha_\eta\eta^{b+1}(1-\zeta)^{-a-b-2}$&\\

$\displaystyle \,\frac{x^ay^bz^c}{(1+z)^{k+2}} \col{x \left( 1-x \right) }{0}{0}$& $\frac{\xi^a\eta^b\zeta^c}{(1-\zeta)^{2+a+b+c-k}  }\col{\xi\alpha_\xi}{0}{0}$ & Volume shape functions,  $0\leq a,c\leq k-2, 0\leq b\leq k-1$.\\

$\displaystyle \,\frac{x^ay^bz^c}{(1+z)^{k+2}} \col{0}{y \left( 1-y\right)}{0}$& $\frac{\xi^a\eta^b\zeta^c}{(1-\zeta)^{2+a+b+c-k}} \col{0}{\eta\alpha_\eta}{0}$ & Volume shape functions, $0\leq b,c\leq k-2, 0\leq a\leq k-1$.\\

$\displaystyle \, \frac{x^ay^bz^c}{(1+z)^{k+2}}\col{0}{0}{z}$& $\frac{\xi^a\eta^b\zeta^{c+1}}{(1-\zeta)^{2+a+b+c-k}} \col{-\xi}{-\eta}{1-\zeta}$ & Volume shape functions, $0\leq a,b,c\leq k-2$.\\

   \end{tabular}
   \caption{Shape functions for 2-forms on a pyramid. Since the approximation space $\D(\ip)$ is invariant under the rotation, $R_\infty:\ip\rightarrow \ip$,  it is only necessary to demonstrate shape functions for a representative base vertex, vertical edge, base edge and vertical face. Then, using  \eqref{rotation} and the subsequent remarks, the inverse pullback of these to the finite pyramid will also be invariant under the rotation $R$. Note that $\alpha_\xi:= (1-\zeta-\xi), \alpha_\eta:=(1-\zeta-\eta)$}
   \label{tab:basisfunctions2}
  \end{table}